\RequirePackage{amsthm}

\documentclass[sn-mathphys-num]{sn-jnl}

\usepackage{graphicx}%
\usepackage{multirow}%
\usepackage{amsmath}%
\usepackage{amssymb}%
\usepackage{amsfonts}%
\usepackage{booktabs}%
\usepackage{algorithm}%
\usepackage{algorithmicx}%
\usepackage{mathtools}%
\usepackage{comment}

\newcommand{\N}{\mathbb{N}}
\newcommand{\R}{\mathbb{R}}
\newcommand{\bx}{\boldsymbol{x}}

\DeclareMathOperator{\spn}{span}
\DeclareMathOperator{\Range}{Rg}
\DeclareMathOperator{\dom}{dom}

\theoremstyle{thmstyleone}%
\newtheorem{teo}{Theorem}
\newtheorem{prop}[teo]{Proposition}%
\newtheorem{coro}[teo]{Corollary}%

\theoremstyle{thmstyletwo}%
\newtheorem*{ass}{Conditions}%
\newtheorem{exmp}{Example}%
\newtheorem{rmk}{Remark}%

\theoremstyle{thmstylethree}%

\begin{document}

\title[Article Title]{Convergence analysis and parameter estimation for the iterated Arnoldi-Tikhonov method}

\author[1]{\fnm{Davide} \sur{Bianchi}}\email{bianchid@mail.sysu.edu.cn}

\author[2]{\fnm{Marco} \sur{Donatelli}}\email{marco.donatelli@uninsubria.it}

\author[2]{\fnm{Davide} \sur{Furchì}}\email{dfurchi@uninsubria.it}

\author[3]{\fnm{Lothar} \sur{Reichel}}\email{reichel@math.kent.edu}

\affil[1]{\orgdiv{School of Mathematics (Zhuhai)}, \orgname{Sun Yat-sen University}, \orgaddress{\city{Zhuhai}, \postcode{519082}, \country{China}}}

\affil[2]{\orgdiv{Dipartimento di
		Scienza e Alta Tecnologia}, \orgname{Universit\`a dell'Insubria}, \orgaddress{\city{Como}, \postcode{22100}, \country{Italy}}}

\affil[3]{\orgdiv{Department of
		Mathematical Sciences}, \orgname{Kent State University}, \orgaddress{\city{Kent}, \postcode{OH 44242}, \country{USA}}}

\abstract{The Arnoldi-Tikhonov method is a well-established regularization technique for solving 
large-scale ill-posed linear inverse problems. This method leverages the Arnoldi 
decomposition to reduce computational complexity by projecting the discretized problem
into a lower-dimensional Krylov subspace, in which it is solved. 
This paper explores the iterated Arnoldi-Tikhonov method, conducting a comprehensive analysis that addresses all approximation errors. Additionally, it introduces a novel strategy for choosing the regularization parameter, leading to more accurate approximate solutions compared to the standard Arnoldi-Tikhonov method.
Moreover, the proposed method demonstrates robustness with respect to the regularization parameter, as confirmed by the numerical results.}


\maketitle

\section{Introduction}
Let $T\colon\mathcal{\mathcal{X}}\to \mathcal{\mathcal{Y}}$ be a bounded linear operator between two separable Hilbert spaces $\mathcal{X}$ and $\mathcal{Y}$, and assume that $T$ is not continuously invertible. We are concerned with the solution of operator equations of the form
\begin{equation}\label{eq}
  Tx=y,\quad x\in\mathcal{\mathcal{X}},\quad y\in\mathcal{Y}.
\end{equation}
This model equation is assumed to be solvable; we denote its unique least-squares solution
of minimal norm by $x^\dagger$. Since $T$ is not continuously invertible, $x^\dagger$ 
might not depend continuously on $y$. Solving~\eqref{eq} therefore is an ill-posed 
problem.

In problems of interest to us, the right-hand side of equation~\eqref{eq} is not 
available; only an error-contaminated approximation $y^{\delta}\in\mathcal{Y}$ of $y$ is
known. For instance, we may determine $y^{\delta}$ from measurements. Assume that 
$y^{\delta}$ satisfies
\begin{equation*}
    \|y-y^{\delta}\|_{\mathcal{Y}}\leq\delta,
\end{equation*}
with a known bound $\delta >0$; here $\|\cdot\|_{\mathcal{Y}}$ denotes the norm in $\mathcal{Y}$.

We are therefore led to the problem of finding a solution to the equation
\begin{equation}\label{eqdelta}
    Tx^{\delta}=y^{\delta},\quad \mbox{where } x^{\delta}\in\mathcal{X},\quad y^{\delta}\in\mathcal{Y}.
\end{equation}
However, the least-squares solution of minimal norm of~\eqref{eqdelta} typically is not a meaningful approximation of the solution $x^\dagger$ of~\eqref{eq} due to the fact that $T$ does not have a  continuous inverse. Given $T$ and $y^\delta$, we would like to determine an accurate 
approximation of $x^\dagger$. This requires that equation~\eqref{eqdelta} be 
\emph{regularized}, i.e., that the operator $T$ be replaced with a nearby operator, so that the computation of the solution of the equation with the nearby operator is a 
well-posed problem. This solution is less sensitive to perturbations in the right-hand
side than the solution of~\eqref{eqdelta}. Equations of the form~\eqref{eqdelta} arise in
a variety of applications including remote sensing~\cite{diazdealba2019}, atmospheric 
tomography~\cite{ramlau2012}, computerized tomography~\cite{natterer2001}, adaptive optics 
\cite{raffetseder2016}, and image restoration~\cite{bentbib2018}.

A significant portion of the literature on the solution of linear operator equations with
an operator with a discontinuous inverse, such as~\eqref{eq} or~\eqref{eqdelta}, primarily 
concentrates on analyzing the properties of these problems in infinite-dimensional Hilbert 
spaces. However, properties of the finite-dimensional problem that emerges through the 
unavoidable discretization of the model equation are frequently overlooked. Conversely, 
many works address solution techniques for the discretized problem, which may be large in 
scale, but neglect to account for the effects of discretization and approximation errors 
when reducing the dimension of the problem.

Based on results by Natterer~\cite{natterer1977} and Neubauer~\cite{neubauer1988a}, Ramlau
and Reichel in~\cite{ramlau2019} addressed the aforementioned gap within the setting of 
Arnoldi-Tikhonov methods 
(\cite{calvetti2000tikhonov,lewis2009arnoldi,gazzola2014generalized,gazzola2015krylov}). 
First the operator equation~\eqref{eqdelta} is discretized, and the effect of the 
discretization is established based on the analysis by Natterer~\cite{natterer1977}. The 
resulting linear system of equations is assumed to be large. The matrix representing this
system is subsequently reduced in size through the application of an Arnoldi decomposition. 
The reduced linear system of equations is regularized using the Tikhonov method and solved
straightforwardly. The error stemming from replacing the original system by a smaller one 
is examined using the results presented in~\cite{neubauer1988a}.

It is the purpose of this paper to extend the analysis in~\cite{natterer1977}, 
\cite{neubauer1988a}, and~\cite{ramlau2019} to the iterated version of the 
Arnoldi-Tikhonov method. It is a well-established fact that iterative Tikhonov methods 
frequently produce computed approximate solutions of superior quality, with a higher 
robustness of the method, when compared to standard Tikhonov regularization; see 
\cite{buccini2017iterated,donatelli2012nondecreasing,hanke1998nonstationary}. A first 
approach to combine an Arnoldi approximation with the iterated Tikhonov method was
described by Buccini et al.~\cite{buccini2023arnoldi}, where an Arnoldi-based 
preconditioner is  employed in a preconditioned iteration method proposed in 
\cite{donatelli2013fast}, and inspired by the iterated Tikhonov method. 

Differently from preconditioning and hybrid Arnoldi-Tikhonov methods, this paper uses the 
Arnoldi approximation to develop the iterated Tikhonov method.  We carry out a 
comprehensive analysis that takes into account all approximation errors. To the best of 
our knowledge, this is the first time such an analysis has been performed within the framework
of iterated Tikhonov regularization. Our analysis leads to a new approach to determining the
regularization parameter and improves on the parameter choice method discussed in 
\cite{ramlau2019}, in the sense that it gives computed approximate solutions of higher 
quality.

This paper is organized as follows. Section~\ref{sec:pre} describes the setting of this 
work. The section discusses the discretization of the model problem and presents the basic
Tikhonov regularization method. The Arnoldi-Tikhonov method is described in 
Section~\ref{sec:AT}, and Section~\ref{sec:iAT} presents the iterated Arnoldi-Tikhonov 
method and provides convergence results. Section~\ref{sec:comp} presents a few computed 
examples, including a comparison of the Arnoldi-Tikhonov method in~\cite{ramlau2019} and 
our iterated Arnoldi-Tikhonov method. Concluding remarks can be found in 
Section~\ref{sec:end}. Appendix~\ref{sec:appA} provides all the technical details of the 
theory that underpin the results of this work.


\section{Preliminary notations and assumptions}\label{sec:pre}
For a rigorous introduction to regularization theory for inverse problems in Hilbert spaces, 
we refer readers to~\cite{engl1996,scherzer2009variational}. Let us fix a bounded linear 
operator $T\colon\mathcal{\mathcal{X}}\to \mathcal{\mathcal{Y}}$ between two separable 
Hilbert spaces $\mathcal{X}$ and $\mathcal{Y}$, and let 
$T^{\ast}\colon\mathcal{Y}\to\mathcal{X}$ denote its adjoint. Let 
$\|\cdot\|_{\mathcal{X}}$ and $\|\cdot\|_{\mathcal{Y}}$ stand for the norms induced by the 
inner products of $\mathcal{X}$ and $\mathcal{Y}$, respectively. When dealing with 
distinct Hilbert spaces, the norm will be denoted with the name of the respective space as
a subscript. Moreover we let $\|T\|$ denote the operator norm of a bounded 
linear operator $T$. In the particular case of the Euclidean norm, we will use the standard 
notation $\|\cdot\|_2$ also for the norm induced on bounded operators.

Let $T^\dagger$ stand for the Moore-Penrose pseudo-inverse of $T$, that is
\begin{equation*}
	T^\dagger \colon \dom(T^\dagger)\subseteq \mathcal{Y} \to \mathcal{X}, \qquad \mbox{where } \dom(T^\dagger)= \Range(T)\oplus \Range(T)^\perp.
\end{equation*}
For any $y\in\dom(T^\dagger)$, the element $x^\dagger\coloneqq T^\dagger y$ is the unique
least-square solution of minimal norm of the model equation~\eqref{eq}; it is referred to 
as the \emph{best-approximate solution}. To ensure consistency in~\eqref{eq}, we will 
assume as a base hypothesis that
\begin{equation*}
	y\in \dom(T^\dagger).
\end{equation*}

Since $T$ is not continuously invertible, the operator $T^\dagger$ is unbounded. This may make the
least-squares solution $T^\dagger y$ of~\eqref{eq} very sensitive to perturbations in
$y$. A regularization method replaces $T^\dagger$ by a member of a family 
$\{R_\alpha \colon \mathcal{Y} \to \mathcal{X}\}$ of continuous operators that depend on a
parameter $\alpha$, paired with a suitable parameter choice rule 
$\alpha=\alpha(\delta, y^\delta)>0$. Roughly, the pair $(R_\alpha,\alpha)$ furnishes a 
point-wise approximation of $T^\dagger$; see~\cite[Definition 3.1]{engl1996} for a rigorous
definition.

\subsection{Discretization and Tikhonov regularization}\label{ssec:disc}
In real-world applications, the model equations~\eqref{eq}-\eqref{eqdelta} 
are discretized before computing an approximate solution. The discretization process 
introduces a discretization error. To bound the propagated discretization error, we use 
results from Natterer~\cite{natterer1977} and follow a similar approach as in 
\cite{ramlau2019}, but applied to the iterative version of the Tikhonov method. 

Consider a sequence $\mathcal{X}_1 \subset \mathcal{X}_2 \subset \ldots \subset 
\mathcal{X}_n \subset \ldots \subset \mathcal{X}$ of finite-dimensional subspaces, whose
union is dense in $\mathcal{X}$, with $\dim(\mathcal{X}_n)=n$. Define the projectors 
$P_n\colon\mathcal{\mathcal{X}}\to\mathcal{\mathcal{X}}_n$ and 
$Q_n\colon\mathcal{Y}\to \mathcal{Y}_n \coloneqq T(\mathcal{X}_n)$, and the inclusion 
operator $\iota_n \colon \mathcal{X}_n \hookrightarrow \mathcal{X}$. Application of these 
operators to equations~\eqref{eq} and~\eqref{eqdelta} gives the equations
\begin{align*}
	&Q_nT\iota_n P_n x=Q_ny,\\
    &Q_nT\iota_nP_nx^{\delta}=Q_ny^{\delta}.
\end{align*}
Define the operator $T_n \colon \mathcal{X}_n \to \mathcal{Y}_n$,
$$
T_{n}\coloneqq Q_nT\iota_n,
$$
and the finite-dimensional vectors 
\begin{equation*}
    y_n \coloneqq Q_ny,\quad y^\delta_n \coloneqq Q_ny^\delta,\quad x_n \coloneqq P_nx,\quad x_n^{\delta} \coloneqq P_nx^{\delta}.
\end{equation*}
It is natural to identify $T_{n}$ with a matrix in $\R^{n\times n}$, and $y_n$, 
$y^\delta_n$, $x_n$, and $x_n^{\delta}$ with elements in $\R^n$. This gives us the linear systems of 
equations
\begin{align}
	&T_{n}x_n=y_n,\label{disceq2}\\
    &T_{n}x_n^{\delta}=y_n^{\delta}.\label{disceqdelta2}
\end{align}
From now on we will consider $T_n$ a square matrix, not necessarily injective, that 
represents a discretization of $T$. 

Let $T^{\dagger}_{n}$ denote the Moore-Penrose pseudoinverse of the matrix $T_{n}$. Then 
the unique least-squares solutions with respect to the Euclidean vector norm  of equations
\eqref{disceq2} and~\eqref{disceqdelta2} are given by
\begin{equation*}
  x_n^{\dagger}\coloneqq T^{\dagger}_{n}y_n\quad\mbox{and}\quad 
  x_n^{\dagger, \delta}\coloneqq T^{\dagger}_{n}y_n^{\delta},
\end{equation*}
respectively. The fact that the operator $T$ has an unbounded inverse may result in 
that the matrix $T_n$ is severely ill-conditioned. Consequently, $x_n^{\dagger, \delta}$ may be
far from being an accurate approximation of $x_n^{\dagger}$. It follows that 
regularization of the discretized operator equation~\eqref{disceqdelta2} is required; see 
the end of this section and Section~\ref{sec:AT}. 

Note also that the best approximate solution $x^{\dagger}_n\in\mathcal{X}_n$ of~\eqref{disceqdelta2} might 
not be an accurate approximation of the desired solution $x^\dagger$ of~\eqref{eq}, due to
a large propagated error stemming from the discretization. We therefore would like to 
determine a bound for $\|x^\dagger-x^{\dagger}_n\|_{\mathcal{X}}$. This is in general not
possible without some additional assumptions. In particular, it is not sufficient for $T$
and $T_{n}$ to be just close in the operator norm; see~\cite[Example 3.19]{engl1996}. We 
will therefore assume that  
\begin{equation}\tag{H1}\label{hp:x_n-convergence}
\|x^\dagger-x^{\dagger}_n\|_{\mathcal{X}}\leq f(n)\rightarrow 0\quad\mbox{as } n\to\infty,
\end{equation}
for a suitable function $f$. For example, if $T$ is compact and 
$\limsup_{n\to \infty}\|T_n^{\dagger *}x^\dagger_n\|_{\mathcal{X}}<\infty$, or if $T$ is 
compact and $\{\mathcal{\mathcal{X}}_n\}_n$ is a discretization resulting from the 
\emph{dual-least square projection} method, see~\cite[Section 3.3]{engl1996}, then 
\begin{equation*}
f(n)=O(\|(I-P_n)T^{\ast}\|),
\end{equation*}
where $I$ denotes the identity map.

Another scenario arises when $T$ is injective. In such a case, establishing~\eqref{hp:x_n-convergence} by substituting the norm of $\mathcal{X}$ with the norm induced by $T$  directly ensures the robustness of a general projection method. This approach is detailed in~\cite[Satz 2.2]{natterer1977}. For practical applications and specific asymptotic bounds related to this case, see
\cite[Section 4]{natterer1977} and~\cite[Section 2]{ramlau2019}.

Finally, let $\lbrace e_j\rbrace_{j=1}^n$ form a convenient basis for $\mathcal{X}_n$ and
consider the representation
\begin{equation*}
  x_n=\sum\limits_{j=1}^nx_j^{(n)}e_j
\end{equation*}
of an element $x_n\in\mathcal{\mathcal{X}}_n$. We identify the element $x_n$ with the 
vector 
\begin{equation*}
  \bx_n=[x_1^{(n)},\ldots,x_n^{(n)}]^{\ast}\in\R^n,
\end{equation*}
where the superscript $\ast$ indicates transposition for matrices. If $\lbrace e_j\rbrace_{j=1}^n$ is an orthonormal basis, then we may choose the norms so
that $\|\bx_n\|_2 = \|x_n\|_{\mathcal{X}}$. However, for certain discretization methods, 
the basis is not orthonormal and this equality does not hold. Hence, as in 
\cite{ramlau2019}, we make the assumption that there exist positive constants $c_{\min}$
and $c_{\max}$, independent of $n$, such that
\begin{equation}\tag{H2}\label{hp:norm_equivalence}
	c_{\min}\|\bx_n\|_2\leq\|x_n\|_\mathcal{X}\leq c_{\max}\|\bx_n\|_2.
\end{equation}
This condition holds in many practical scenarios, such as wavelets, the discrete cosine transform~\cite{goodman2016discrete}, and linear B-splines ~\cite{de1978practical}. In particular, $c_{\min} = c_{\max}  = 1$ if $\lbrace e_j\rbrace_{j=1}^n$ forms  an orthonormal basis; see~\cite{ramlau2019} for practical examples and further details.

Let $\alpha>0$ be a regularization parameter. Tikhonov regularization in standard form 
applied to~\eqref{disceqdelta2} reads 
\begin{equation*}
x_{\alpha,n}^{\delta}\coloneqq \underset{x_n \in \R^n}{\operatorname{argmin}} 
\|T_nx_n-y_{n}^{\delta}\|_2^2+\alpha\|x_{n}\|_2^2,
\end{equation*}
whose closed form solution is given by 
\begin{equation}\label{Tikhonov_n:closed_form}
	x_{\alpha,n}^{\delta}= (T_n^{\ast}T_n+\alpha I_n)^{-1}T_n^{\ast}y_n^\delta,
\end{equation}
where $I_n$ denotes the identity matrix of order $n$. Clearly, when $n$ is large, computing $x_{\alpha,n}^{\delta}$ by
using formula~\eqref{Tikhonov_n:closed_form} is impractical. In the next section, we will
discuss how to reduce the complexity of Tikhonov regularization and achieve a fairly 
accurate approximation of $x_{\alpha,n}^{\delta}$.


\section{The Arnoldi-Tikhonov method}\label{sec:AT}
When working with large square matrices, the Arnoldi decomposition is a commonly used 
technique to reduce the computational complexity, while retaining some of the important 
characteristics of these matrices. For a comprehensive review of the Arnoldi process, we
refer to~\cite[Section 6.3]{saad2003}; an algorithm that implements the Arnoldi process 
for computing the Arnoldi decomposition also can be found in~\cite{ramlau2019}. This 
section reviews the Arnoldi decomposition and shows how it can be applied to approximate 
the discretized equation~\eqref{disceqdelta2} and its Tikhonov regularized solution 
\eqref{Tikhonov_n:closed_form} in a low-dimensional subspace $\R^\ell$, with 
$1\leq \ell\ll n$.

\subsection{The Arnoldi approximation}
Application of $1\leq \ell\ll n$ steps of the Arnoldi process to the matrix $T_{n}$, with
initial vector $y_n^{\delta}$, gives the Arnoldi decomposition
\begin{equation}
	T_{n}V_{n,\ell}=V_{n,\ell+1}H_{\ell+1,\ell}.\label{arn}
\end{equation}
The columns of the matrix
\begin{equation*}
V_{n,\ell+1}=\left[\begin{array}{@{}c|c|c@{}}
	v_{n,\ell+1}^{(1)} &\cdots & v_{n,\ell+1}^{(\ell+1)}
\end{array}\right]=
\left[\begin{array}{@{}c|c@{}}
	V_{n,\ell} & v_{n,\ell+1}^{(\ell+1)}
\end{array}\right]\in\R^{n\times (\ell+1)}
\end{equation*}
form an orthonormal basis for the Krylov subspace
\begin{equation*}
\mathcal{K}_{\ell+1}(T_{n},y_n^{\delta})=
\spn\lbrace y_n^{\delta},T_{n}y_n^{\delta},\ldots,T_{n}^{\ell}y_n^{\delta}\rbrace
\end{equation*}
with respect to the canonical inner product, that is, 
$V_{n,\ell+1}^{\ast}V_{n,\ell+1} = I_{\ell+1}$. Moreover, 
$H_{\ell+1,\ell}\in\R^{(\ell+1)\times\ell}$ is an upper Hessenberg matrix, i.e., all
entries below the subdiagonal vanish. Clearly, both $V_{n,\ell+1}$ and  
$H_{\ell+1,\ell}$ depend on both $T_n$ and $y_n^\delta$.

If the Arnoldi process does not break down, then the matrix $H_{\ell+1,\ell}$ has full rank. In rare situations, the Arnoldi process breaks down at step $j\leq\ell$. If this 
happens, then the decomposition~\eqref{arn} reduces to 
\begin{equation*}
	T_{n}V_{n,j}=V_{n,j}H_{j,j},
\end{equation*}
and the solution of~\eqref{disceqdelta2} lives in the Krylov subspace 
$\mathcal{K}_{j}(T_{n},y_n^{\delta})$ if $H_{j,j}$ is nonsingular, which is guaranteed 
when $T_n$ is nonsingular. There exist techniques to continue the Arnoldi process when it
breaks down and $H_{j,j}$ is singular, see for example~\cite{reichel2005breakdown}, but 
this is out of the scope of this paper. Henceforth, we shall proceed with the understanding 
that the Arnoldi decomposition~\eqref{arn} is at our disposal.

Define the following approximation of the matrix $T_n$,
\begin{equation}\label{T_n-approximation}
 T_n^{(\ell)}\coloneqq V_{n,\ell+1}H_{\ell+1,\ell}V_{n,\ell}^{\ast} \in \R^{n\times n},
\end{equation}
which we will refer to as the \emph{Arnoldi approximation} of $T_n$. From now on, we shall
operate under the assumption that
\begin{equation}\label{hell}
	\|T_n-T_n^{(\ell)}\|_2\leq h_{\ell} \quad \mbox{for some} \quad h_{\ell}\geq0.
\end{equation}
Notice that $T_n^{(\ell)} = T_nV_{n,\ell}V_{n,\ell}^{\ast}$, which reduces to $T_n$ when 
$\ell=n$, since in this case $V_{n,\ell}V_{n,\ell}^{\ast}=I_n$. 

\subsection{The Arnoldi-Tikhonov method}
Having introduced the Arnoldi approximation $T_n^{(\ell)}$ in~\eqref{T_n-approximation}, 
we move from the discretized equation~\eqref{disceqdelta2} to the \emph{approximated} 
discretized equation
\begin{equation}\label{eq:approximated_discretized_equation}
  T_n^{(\ell)}x_n^{\delta}= y_n^\delta,  
\end{equation}
and its associated Tikhonov regularized solution 
\begin{equation}\label{arntm}
x_{\alpha,n}^{\delta,\ell}\coloneqq (T_n^{(\ell)*}T_n^{(\ell)}+\alpha I_n)^{-1}
T_n^{(\ell)*}y_n^{\delta}.
\end{equation}
We can reduce the computational complexity of solving~\eqref{arntm} by exploiting the 
structure of $T_n^{(\ell)}$. Observe that, for any $x_\ell \in \R^\ell$,
\begin{align*}
V_{n,\ell} (H^{\ast}_{\ell+1,\ell}H_{\ell+1,\ell}+\alpha I_{\ell})x_\ell &= (V_{n,\ell}H^{\ast}_{\ell+1,\ell}V^{\ast}_{n,\ell+1}V_{n,\ell+1}H_{\ell+1,\ell}
V^{\ast}_{n,\ell}+\alpha I_{n})V_{n,\ell}x_\ell\\
&= (T_{n}^{(\ell)*}T_{n}^{(\ell)}+\alpha I_n)V_{n,\ell}x_\ell
\end{align*}
and, hence, it holds
\begin{equation}\label{HV}
    (T_{n}^{(\ell)*}T_{n}^{(\ell)}+\alpha I_n)^{-1}V_{n,\ell}=V_{n,\ell}(H^{\ast}_{\ell+1,\ell}H_{\ell+1,\ell}+\alpha I_{\ell})^{-1}.
\end{equation}
Therefore, if we define 
$y_{\ell+1}^\delta\coloneqq V_{n,\ell+1}^{\ast}y_n^\delta\in\R^{\ell+1}$ and
\begin{equation*}
z_{\alpha,\ell}^{\delta,\ell}\coloneqq (H^{\ast}_{\ell+1,\ell}H_{\ell+1,\ell}+
\alpha I_{\ell})^{-1}H^{\ast}_{\ell+1,\ell}y_{\ell+1}^\delta, 
\end{equation*}
which is the Tikhonov regularized solution associated with the \emph{reduced} approximated
discretized equation
\begin{equation}\label{eq:projected_approximated_discretized_equation}
    H_{\ell+1,\ell}z_\ell = y_{\ell+1}^\delta.
\end{equation}
Combining equations~\eqref{T_n-approximation},~\eqref{arntm} and~\eqref{HV} we obtain
\begin{equation*}
   x_{\alpha,n}^{\delta,\ell} =  V_{n,\ell}z_{\alpha,\ell}^{\delta,\ell}.
\end{equation*}
That is, the Tikhonov regularized solution $x_{\alpha,n}^{\delta,\ell}$ of equation 
\eqref{eq:approximated_discretized_equation} is the back-projection in $\R^n$, through 
$V_{n,\ell}$, of the Tikhonov regularized solution 
$z_{\alpha,\ell}^{\delta,\ell}\in \R^\ell$ 
of the low-dimensional reduced equation 
\eqref{eq:projected_approximated_discretized_equation}. This procedure of computing an 
approximate solution of~\eqref{disceqdelta2}, and therefore of~\eqref{eqdelta}, is 
referred to as the \emph{Arnoldi-Tikhonov} (AT) \emph{regularization method}. We summarize 
it in the following algorithm. 

\begin{algorithm}
\caption{The Arnoldi-Tikhonov regularization method}\label{algo:AT}
\begin{algorithmic}[1]
\State \textbf{Input:} $\{T_n$, $y^\delta_n$, $\ell\}$
\State \textbf{Output:} $x_{\alpha,n}^{\delta,\ell}$
\State Compute $\{V_{n,\ell+1}$, $H_{\ell+1,\ell}\}$ with the Arnoldi decomposition 
algorithm~\cite[Section 6.3]{saad2003}
\State Compute $y^\delta_{\ell+1}= V^{\ast}_{n,\ell+1}y^\delta_n$ 
\State Determine $\alpha$
\State Compute $z_{\alpha,\ell}^{\delta,\ell}= (H^{\ast}_{\ell+1,\ell}H_{\ell+1,\ell}+
\alpha I_{\ell})^{-1}H^{\ast}_{\ell+1,\ell}y_{\ell+1}^\delta$
\State Return $x_{\alpha,n}^{\delta,\ell}=V_{n,\ell}z_{\alpha,\ell}^{\delta,\ell}$
\end{algorithmic}
\end{algorithm}

Note that step $6$ in Algorithm \ref{algo:AT} does not require computing the inverse of the matrix $H^{\ast}_{\ell+1,\ell}H_{\ell+1,\ell}+
\alpha I_{\ell}$. Instead, $z_{\alpha,\ell}^{\delta,\ell}$ can be computed by using the 
Cholesky factorization of  $H^{\ast}_{\ell+1,\ell}H_{\ell+1,\ell}+ \alpha I_{\ell}$, 
which can be updated from the factorization computed at the previous step $\ell-1$.

A convergence analysis for approximate solutions of~\eqref{eqdelta} computed with Algorithm~\ref{algo:AT} is carried out in~\cite{ramlau2019}, where also convergence rates for 
$\|x_n^\dagger - x_{\alpha,n}^{\delta,\ell}\|_2$ are established. Convergence in the space
$\mathcal{X}$ then is obtained by using~\eqref{hp:norm_equivalence}.


\section{The iterated Arnoldi-Tikhonov (iAT) method}\label{sec:iAT}
This section discusses an iterative extension of the AT regularization method
described by 
Algorithm \ref{algo:AT}. We refer to this extension as the iterated Arnoldi-Tikhonov (iAT)
method and present an analysis of its convergence properties.

Our motivation for introducing an iterated version of the AT method is twofold. Firstly, 
the standard Tikhonov regularization method exhibits a saturation rate of 
$O(\delta^{2/3})$ as $\delta\searrow 0$, as demonstrated in 
\cite[Proposition 5.3]{engl1996}. This implies that approximate solutions computed with 
the AT method do not converge to the solution of~\eqref{eqdelta} faster than 
$O(\delta^{2/3})$ as $\delta\searrow 0$, a behavior that also is suggested by 
\cite[Corollary 4.6]{ramlau2019}. However, the iAT method can surpass this saturation 
rate, as we demonstrate in Corollary \ref{convtot} below.

Secondly, the iterated Tikhonov method typically yields significantly improved approximate 
solutions, compared with approximate solutions determined by the noniterative AT method.
This is, in particular, the case when iAT is combined with the parameter choice rule of 
Proposition \ref{convf} below. Computed illustrations are presented in Section 
\ref{sec:comp}. We remark that the improved quality of the computed solutions determined
by iAT is achieved for essentially the same computational cost as required by the
noniterative AT method.

It is worth noting that the parameter choice rule $\alpha=\alpha(\delta, y_n^\delta)$ 
introduced in Proposition \ref{convf} enhances the quality of computed solutions also when 
applied with the noniterative AT method. The parameter choice rule of Proposition 
\ref{convf} is optimal in the sense that the quality of the computed solution does not 
improve by choosing a regularization parameter larger than $\alpha$. This is shown 
in Proposition \ref{comp_iAT} below.

\subsection{Definition of the iAT method}
We combine iterated Tikhonov regularization with the AT method. Given the operator 
$T_n^{(\ell)}$ in~\eqref{T_n-approximation}, we define the iteration 
\begin{equation}\label{iATiter}
  \begin{cases}
    (T_n^{(\ell)\ast}T_n^{(\ell)}+\alpha I_n)x_{\alpha,n,i}^{\delta,\ell}=
    T_{n}^{(\ell)\ast}y_n^{\delta}+\alpha x_{\alpha,n,i-1}^{\delta,\ell},
        \qquad i \geq 1,\\
    x_{\alpha,n,0}^{\delta,\ell}\coloneqq 0,
  \end{cases}
\end{equation}
which in variational form reads
\begin{equation*}
  \begin{cases}
    x_{\alpha,n,i}^{\delta,\ell}\coloneqq \underset{x_n \in \R^n}{\operatorname{argmin}}
    \|T_n^{(\ell)}x_n-y_n^\delta \|_2^2+\alpha\|x_n-x_{\alpha,n,i-1}^{\delta,\ell}\|_2^2,
        \qquad i \geq 1,\\
    x_{\alpha,n,0}^{\delta,\ell}\coloneqq 0.
  \end{cases}
\end{equation*}
In particular, 
\begin{equation*}
x_{\alpha,n,i}^{\delta,\ell}=\sum\limits_{k=1}^i\alpha^{k-1}(T_{n}^{(\ell)\ast}
T_n^{(\ell)}+\alpha I_n)^{-k}T_{n}^{(\ell)\ast}y_n^{\delta}.
\end{equation*}
We will denote the analogous solution when the vector $y_n$ replaces $y_n^{\delta}$ 
in the above equation by $x_{\alpha,n,i}^{\ell}$. Similarly to the discussion in
Section \ref{sec:AT}, we can leverage the Arnoldi decomposition to establish that
\begin{equation*}
  x_{\alpha,n,i}^{\delta,\ell}=V_{n,\ell}z_{\alpha,\ell,i}^{\delta,\ell},
\end{equation*}
where
\begin{equation*}
  z_{\alpha,\ell,i}^{\delta,\ell}\coloneqq\sum_{k=1}^i\alpha^{k-1}(H_{\ell+1,\ell}^{\ast}
  H_{\ell+1,\ell}+\alpha I_{\ell})^{-k}H_{\ell+1,\ell}^{\ast}y_{\ell+1}^{\delta}.
\end{equation*}
We summarize the iAT method in Algorithm \ref{algo:itarntik} below. Step 6 of the 
algorithm is evaluated by using an iteration analogous to~\eqref{iATiter} and applying a 
Cholesky factorization of the matrix 
$H_{\ell+1,\ell}^{\ast} H_{\ell+1,\ell}+\alpha I_{\ell}$. Notice that for $i=1$, we 
recover the AT method of Algorithm \ref{algo:AT}. When the matrix $T_n$ is large and $\ell\ll n$, which we
assume to be the case, the main computational effort required by Algorithm 
\ref{algo:itarntik} is the evaluation of the $\ell$ matrix-vector products with the matrix 
$T_n$ required to compute the Arnoldi decomposition~\eqref{arn}. In particular, the 
computational effort required by Algorithm \ref{algo:itarntik} is essentially independent 
of the number of iterations $i$ of the iAT method.

\begin{algorithm}
\caption{The iterated Arnoldi-Tikhonov method}\label{algo:itarntik}
\begin{algorithmic}[1]
\State \textbf{Input:} $\{T_n$, $y^\delta_n$, $\ell$, $i\}$
\State \textbf{Output:} $x_{\alpha,n,i}^{\delta,\ell}$
\State Compute $\{V_{n,\ell+1}$, $H_{\ell+1,\ell}\}$ with the Arnoldi process 
\cite[Section 6.3]{saad2003}
\State Compute $y^\delta_{\ell+1}= V^{\ast}_{n,\ell+1}y^\delta_n$ 
\State Determine $\alpha$
\State Compute $z_{\alpha,\ell,i}^{\delta,\ell}=\sum\limits_{k=1}^i\alpha^{k-1}
(H_{\ell+1,\ell}^{\ast}H_{\ell+1,\ell}+\alpha I_{\ell})^{-k}H_{\ell+1,\ell}^{\ast}
y_{\ell+1}^{\delta}$
\State Return $x_{\alpha,n,i}^{\delta,\ell}=V_{n,\ell}z_{\alpha,\ell,i}^{\delta,\ell}$
\end{algorithmic}
\end{algorithm}

\subsection{Convergence results}
This section collects convergence results for the iAT method described by Algorithm ~\ref{algo:itarntik}. 
To make this section more readable, we have moved the technical 
details to Appendix \ref{sec:appA}. This section is an adaption of results in 
\cite[Section 3]{ramlau2019}, which, in turn, are based on work by Neubauer
\cite{neubauer1988a}. The connection with the latter paper is made clear in Appendix~\ref{sec:appA}.

We will need the orthogonal projector $\mathcal{R}_{\ell}$ from $\R^n$ into 
$\Range(T_n^{(\ell)})$. This operator also was used in~\cite[Proposition 4.7]{ramlau2019}.
Let $q = \operatorname{rank}(H_{\ell+1,\ell})$ and introduce the singular value 
decomposition 
\begin{equation*}
    H_{\ell+1,\ell}=U_{\ell+1}\Sigma_{\ell+1,\ell}S_{\ell}^{\ast},
\end{equation*} 
where the matrices $U_{\ell+1}\in\R^{(\ell+1)\times(\ell+1)}$ and 
$S_\ell\in\R^{\ell\times\ell}$ are orthogonal, and the nontrivial entries of the diagonal
matrix
\begin{equation*}
    \Sigma_{\ell+1,\ell}={\rm diag}[\sigma_1,\sigma_2,\ldots,\sigma_\ell]\in
\R^{(\ell+1)\times\ell}
\end{equation*}
are ordered according to 
$\sigma_1\geq\ldots\geq\sigma_q>\sigma_{q+1}=\ldots=\sigma_\ell=0$. We note that, since we assume that the Arnoldi process does not break 
down, the matrix $H_{\ell+1,\ell}$ has full rank $q=\ell$. Let
\begin{equation*}
    I_{\ell,\ell+1}=\begin{bmatrix}
  I_\ell & 0\\
  0& 0
  \end{bmatrix}\in\R^{(\ell+1)\times(\ell+1)}.
\end{equation*}
Then, from 
\cite[Proposition 4.7]{ramlau2019}, it holds
\begin{equation*}
    \mathcal{R}_\ell=V_{n,\ell+1}U_{\ell+1}I_{\ell,\ell+1}U_{\ell+1}^*V_{n,\ell+1}^*.
\end{equation*}

Define 
\begin{equation*}
    \hat{y}_{\ell+1}^\delta \coloneqq I_{\ell,\ell+1}U_{\ell+1}^{\ast}y_{\ell+1}^\delta
\end{equation*}
and assume that at least one of the first $q=\ell$ entries of the vector 
$\hat{y}_{\ell+1}^\delta$ is nonvanishing. Then the equation
\begin{equation}\label{condn}
\alpha^{2i+1}(\hat{y}_{\ell+1}^\delta)^{\ast}(\Sigma_{\ell+1,\ell}
\Sigma_{\ell+1,\ell}^{\ast}+\alpha I_{\ell+1})^{-2i-1}\hat{y}_{\ell+1}^\delta=
(Eh_{\ell}+C\delta)^2
\end{equation}
has a unique solution $\alpha>0$. Here $h_{\ell}$ satisfies~\eqref{hell}, and $C$ and $E$ are positive constants such that
\begin{equation}\label{condEC2}
 0\leq Eh_{\ell}+C\delta\leq\|\mathcal{R}_{\ell}y_n^{\delta}\|_2=
 \|U_{\ell+1}I_{\ell,\ell+1}U_{\ell+1}^{\ast}y_{\ell+1}^{\delta}\|_2.
\end{equation}
Equation~\eqref{condn}, which follows from equation~\eqref{select} in Appendix~\ref{sec:appA},
is obtained similarly as~\cite[Proposition 4.8]{ramlau2019}.

\begin{prop}\label{comp_iAT}
Let $C=1$ and $E=\|x_n^{\dagger}\|_2$ and assume that~\eqref{condEC2} holds. Let 
$\alpha>0$ be the 
unique solution of~\eqref{condn}. Then for all $\tilde{\alpha}\geq\alpha$, we have that 
$\|x_n^{\dagger}-x_{\alpha,n,i}^{\delta,\ell}\|_2\leq\|x_n^{\dagger}-
x_{\tilde{\alpha},n,i}^{\delta,\ell}\|_2$.
\end{prop}

\begin{proof}
The result follows from Proposition~\ref{comp} in Appendix~\ref{sec:appA}. 
\end{proof}

\begin{rmk}
Since the left-hand side of~\eqref{condn} is monotonically increasing with
$\alpha>0$, it follows from the 
choice of $E$ in Proposition~\ref{comp_iAT} that the parameter choice strategy 
\eqref{condn} applied with $i=1$, i.e., to Algorithm~\ref{algo:AT}, improves the one 
presented in~\cite{ramlau2019}. This is confirmed by the numerical results in 
Section~\ref{sec:comp}.
\end{rmk}

\begin{prop}\label{convf}
Let $C=1$ and $E=\|x^{\dagger}_n\|_2$. Assume that~\eqref{condEC2} holds and let 
$\alpha>0$ be the unique solution of~\eqref{condn}. Moreover, for some $\nu\geq 0$ and
$\rho>0$, let $x_n^{\dagger}\in\mathcal{X}_{n,\nu,\rho}$, where 
\begin{equation*}
\mathcal{X}_{n,\nu,\rho}\coloneqq\lbrace x_n\in\mathcal{X}_n \mid 
x_n=(T_n^{\ast}T_n)^{\nu}w_n,\; w_n\in \ker(T_n)^\perp \mbox{ and } 
\|w_n\|_2\leq\rho\rbrace.
\end{equation*}
Then
\begin{equation*}
   \|x_n^{\dagger}-x_{\alpha,n,i}^{\delta,\ell}\|_2=
    \begin{cases}
    o(1) &\text{if}\quad\nu=0,\\
    o((h_{\ell}+\delta)^{\frac{2\nu i}{2\nu i+1}})+O(\gamma_{\ell}^{2\nu}\|w_n\|_2)&\text{if}\quad 0<\nu<1,\\
    O((h_{\ell}+\delta)^{\frac{2i}{2i+1}})+O(\gamma_{\ell}\|(I_n-\mathcal{R}_{\ell})T_nw_n\|_2) &\text{if}\quad\nu=1,
 \end{cases}
\end{equation*}
where $\gamma_{\ell}\coloneqq \|(I_n-\mathcal{R}_{\ell})T_n\|_2$.
\end{prop}

\begin{proof}
We first verify that Conditions~\ref{conditions} in Appendix~\ref{sec:appA} are 
satisfied. This follows 
from~\cite[Proposition 4.1]{ramlau2019}. Proposition~\ref{comp_iAT} assures that our 
choice of the regularization parameter is the best in the set $[\alpha,+\infty)$. A bound 
for the convergence rate follows from Theorem~\ref{conv2} since hypothesis \eqref{condplus} is trivially satisfied in a finite dimensional setting.
\end{proof}

As discussed at the end of Appendix~\ref{parasele}, if an estimate of 
$\|x_n^{\dagger}\|_2$ is not available, then we may substitute $E$ by the expression 
$D\|x_{\alpha,n,i}^{\delta,\ell}\|$, with a constant $D\geq 1$. With this choice, for 
$\alpha$ satisfying~\eqref{condEC2}, we achieve the same convergence rates as in 
Proposition~\ref{convf}.
    
\begin{coro}\label{convtot}
Assume that $x_n^{\dagger}\in\mathcal{X}_{n,1,\rho}$ and let $\alpha>0$ be the solution of
\eqref{condn}. Then, for $\ell$ such that $h_{\ell}\sim\delta$, we have
\begin{eqnarray}
\label{equiv1}
  \|x_n^{\dagger}-x^{\delta,\ell}_{\alpha,n,i}\|_2&=&O(\delta^{\frac{2i}{2i+1}})\qquad
  \mbox{as }\delta\rightarrow 0, \\
\label{equiv2}
  \|x^{\dagger}-x_{\alpha,n,i}^{\delta,\ell}\|_{\mathcal{X}}&\leq& f(n)+
  O(\delta^{\frac{2i}{2i+1}})\qquad\mbox{as }\delta\rightarrow 0.
\end{eqnarray}
\end{coro}

\begin{proof}
Equation~\eqref{equiv1} follows from Proposition~\ref{convf} and the fact that
\begin{equation*}
   \|(I_n-\mathcal{R}_{\ell})T_n\|_2=\|(I_n-\mathcal{R}_{\ell})(T_n-T_n^{(\ell)})\|_2=
O(\delta).
\end{equation*}
Inequality~\eqref{equiv2} is obtained from the hypotheses~\eqref{hp:x_n-convergence} and 
\eqref{hp:norm_equivalence}.
\end{proof}

\begin{rmk}
Equation~\eqref{equiv1} shows the expected improvement of the convergence rate $O(\delta^{2/3})$ for
standard Tikhonov regularization.
\end{rmk}


\section{Computed examples}\label{sec:comp}

We apply the iAT regularization method to solve several ill-posed operator equations. 
Examples~\ref{ex1} and \ref{ex2} have previously been considered in~\cite{ramlau2019}, 
while Example \ref{ex3} examines a 2D deblurring problem. The numerical results compare 
Algorithm \ref{algo:AT} (the AT method) and Algorithm \ref{algo:itarntik} (the iAT 
method). All computations were carried out using MATLAB R2024a and 
arithmetic with about $15$ significant decimal digits.

The matrix $T_n$ takes on one of two forms: It either represents a discretization of an 
integral operator (as in Examples \ref{ex1} and \ref{ex2}) or serves as a model for a 
blurring operator (as in Example \ref{ex3}). The latter matrix also may be considered a 
discretized integral operator. With a slightly abuse of notation, we let
within this section $x^{\dagger}_n\in\R^n$ be the vector that is a discretization of 
the exact solution of~\eqref{eq}. Its image $y_n=T_nx_n^{\dagger}$ is 
presumed impractical to measure directly. Instead, we know an observable, 
noise-contaminated vector, $y_n^{\delta}\in\R^n$, that is obtained by adding a vector that models noise to 
$y_n$. Let the vector $e_n\in\R^n$ have 
normally distributed random entries with zero mean. We scale this vector
\begin{equation*}
  \hat{e}_n \coloneqq e_n\frac{\xi\|y_n\|_2}{\|e_n\|_2}
\end{equation*}
to ensure a prescribed noise level $\xi>0$. Then we define 
\begin{equation*}
    y_n^\delta \coloneqq y_n + \hat{e}_n.
\end{equation*}
Clearly,
$\delta \coloneqq \| y_n^\delta - y_n\|_2 = \xi\|y_n\|_2$. We fix the value $\xi$ for each
example such that $\delta$ will correspond to $(100\cdot\xi)\%$ of the norm of $y$. To achieve replicability of the numerical examples, we define the ``noise''
deterministically by setting \texttt{seed}=11 in the MATLAB function \texttt{randn()}, 
which generates normally distributed pseudorandom numbers and we use to determine the entries of the vector $e_n$.

The low-rank approximation $T_n^{(\ell)}$ of $T_n$ is computed by the applying $\ell$
steps of the Arnoldi process to the matrix $T_n$ with initial vector 
$v_1=y_n^{\delta}/\|y_n^{\delta}\|_2$. Thus, we first evaluate the Arnoldi decomposition
\eqref{arn} and then define the matrix $T_n^{(\ell)}$ by~\eqref{T_n-approximation}. Note
that this matrix is not explicitly formed.

We determine the parameter $\alpha$ for 
Algorithm \ref{algo:itarntik} by solving equation~\eqref{condn} with $C=1$ and $E=\|x_n^{\dagger}\|_2$, as suggested by Proposition~\ref{comp}. Inequality~\eqref{condEC2} holds for all examples of this section. In other words, $\alpha$ is the unique solution of
\begin{equation}\label{alpha_iAT}
\alpha^{2i+1}(\hat{y}_{\ell+1}^\delta)^{\ast}(\Sigma_{\ell+1,\ell}
\Sigma_{\ell+1,\ell}^{\ast}+\alpha I_{\ell+1})^{-2i-1}\hat{y}_{\ell+1}^\delta=
(\|x_n^{\dagger}\|_2h_{\ell}+\delta)^2.
\end{equation}
The parameter $\alpha$ for 
Algorithm \ref{algo:AT} is chosen as in~\cite{ramlau2019}, i.e., by solving~\eqref{condn}
with $C=1$, $E=3\|x_n^{\dagger}\|_2$, and $i=1$. In view of Proposition \ref{convf}, our 
parameter choice for $i=1$ improves the one used in~\cite{ramlau2019}. Therefore, it is also 
interesting to evaluate the first iteration of the iAT method.

\begin{table}
\caption{Example~\ref{ex1} - Relative error in approximate solutions computed by iAT and $\alpha$ determined by solving $\eqref{alpha_iAT}$ for different
values of $\ell$, with $n=1000$ and $\xi=0.01$, i.e. $\delta = 1\%$. The AT method is applied 
with the parameter $\alpha$ determined as in~\cite{ramlau2019}.}\label{tab1}
\begin{tabular}{cccccc}
 \toprule%
 & \multicolumn{3}{c}{iAT} & & AT \\
 \cmidrule{2-4}
 \cmidrule{6-6}
 $\ell$ & $i$ & $\alpha$ & $\|x_n^{\dagger}-x_{\alpha,n,i}^{\delta,\ell}\|_2/
 \|x_n^{\dagger}\|_2$ &  & $\|x_n^{\dagger}-x_{\alpha,n}^{\delta,\ell}\|_2/\|x_n^{\dagger}\|_2~~(\alpha)$\\
 \midrule
 \multirow{5}*{10}  & 1 & $3.72\cdot 10^0$ & $1.91\cdot 10^{-1}$ & & 
 \multirow{5}*{$3.32\cdot 10^{-1}\ (9.59)$}\\
 & 50 & $2.80\cdot 10^2$ & $1.46\cdot 10^{-1}$ &\\
 & 100 & $3.33\cdot 10^1$ & $2.70\cdot 10^{-2}$ &\\
 & 150 & $1.04\cdot 10^1$ & $2.06\cdot 10^{-2}$ &\\
 & 200 & $5.80\cdot 10^0$ & $1.72\cdot 10^{-2}$ &\\
 \midrule
 \multirow{5}*{20}  & 1 & $2.25\cdot 10^0$ & $1.41\cdot 10^{-1}$ & & 
 \multirow{5}*{$2.28\cdot 10^{-1}\ (4.98)$}\\
 & 50 & $1.80\cdot 10^2$ & $1.08\cdot 10^{-1}$ &\\
 & 100 & $3.33\cdot 10^1$ & $2.69\cdot 10^{-2}$ &\\
 & 150 & $1.04\cdot 10^1$ & $2.06\cdot 10^{-2}$ &\\
 & 200 & $5.80\cdot 10^0$ & $1.77\cdot 10^{-2}$ &\\
 \midrule
 \multirow{5}*{30} & 1 & $2.24\cdot 10^0$ & $1.41\cdot 10^{-1}$ & & 
 \multirow{5}*{$2.28\cdot 10^{-1}\ (4.96)$}\\
 & 50 & $1.80\cdot 10^2$ & $1.08\cdot 10^{-1}$ &\\
 & 100 & $3.33\cdot 10^1$  & $2.69\cdot 10^{-2}$ &\\
 & 150 & $1.04\cdot 10^1$ & $2.06\cdot 10^{-2}$ &\\
 & 200 & $5.80\cdot 10^0$ & $1.77\cdot 10^{-2}$ &\\
 \botrule
\end{tabular}
\end{table}

\begin{figure}[ht]
\centering
\includegraphics[scale=0.35]{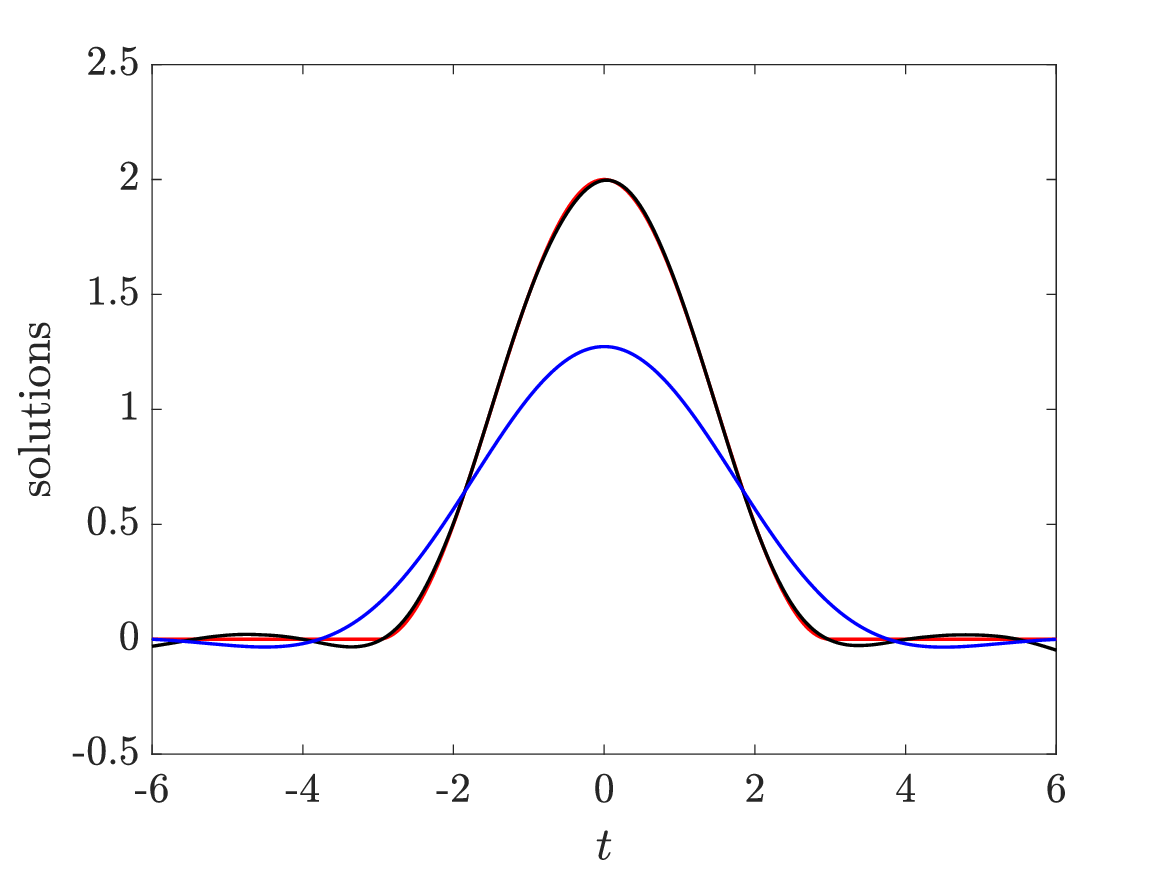}
\caption{Example~\ref{ex1} - Exact solution $x^{\dagger}_n$ (red) and approximate 
solutions $x^{\delta,\ell}_{\alpha,n,i}$ computed by iAT (black) with $i=200$ and $\alpha$ determined by solving~\eqref{alpha_iAT} and 
$x^{\delta,\ell}_{\alpha,n}$ computed by AT (blue) with $\alpha$ determined as in~\cite{ramlau2019}, for $n=1000$, $\ell=10$ and
$\xi=0.01$, i.e. $\delta = 1\%$.}\label{fig1}
\end{figure}

\begin{exmp}\label{ex1}
Consider the Fredholm integral equation of the first kind discussed by Phillips~\cite{phillips1962}:
\begin{equation*}
  \int_{-6}^{6}\kappa(s,t)x(t)dt=y(s),\qquad -6\leq s\leq 6, 
\end{equation*}
for which the solution $x(t)$, the kernel $\kappa(s,t)$, and the right-hand side $y(s)$ are given by
\begin{align*}
 x(t)&=\begin{cases}
        1+\cos(\frac{\pi t}{3}) &\text{if}\quad|t|<3,\\
        0 &\text{if}\quad|t|\geq 3,
       \end{cases}\\
 \kappa(s,t)&=x(s-t),\qquad y(s)=(6-\arrowvert s\arrowvert)\left(1+\frac{1}{2}\cos
 \left(\frac{\pi s}{3}\right)\right)+\frac{9}{2\pi}\sin
 \left(\frac{\pi\arrowvert s\arrowvert}{3}\right).
\end{align*}
We discretize this integral equation by a Nystr\"om method based on the composite 
trapezoidal rule with $n$ nodes using code available in \cite{NRS}. This gives a nonsymmetric matrix $T_n\in\R^{n\times n}$ and the true solution $x_n^{\dagger}\in\R^n$. Table~\ref{tab1} shows the relative error
$\|x_n^{\dagger}-x_{\alpha,n,i}^{\delta,\ell}\|_2/\|x_n^{\dagger}\|_2$ for several
computed approximate solutions $x_{\alpha,n,i}^{\delta,\ell}$ for the noise level $\delta$ of
$1\%$. The error depends on $y_n^{\delta}$, the matrix $T_n$, the approximation error 
\eqref{hell}, and the iteration number $i$. The relative error can be seen to decrease 
slowly as $\ell$ increases for fixed $i$, and to decrease quite rapidly as $i$ increases 
for fixed $\ell$. Moreover, $\alpha$ stabilizes when $i$ increases, independently of
$\ell$. Figure~\ref{fig1} shows the exact solution as well as the approximate solutions 
computed with Algorithms \ref{algo:AT} and \ref{algo:itarntik}. The latter algorithm can 
be seen to determine a much more accurate approximation of the exact solution 
$x_n^{\dagger}$ than the former.

The behavior of the relative error of approximate solutions computed with Algorithm 
\ref{algo:AT} (AT) and Algorithm \ref{algo:itarntik} (iAT) for varying $\alpha$ is 
displayed in Figure~\ref{fig2} in log-log scale. The $\alpha$-values determined by solving
equation~\eqref{condn} are marked by ``$*$'' on the graphs. Note that the regularization
parameters $\alpha$ obtained for the iAT method correspond to a point close to the minimum
of the black curves, while the regularization parameters $\alpha$ for the AT method do not
correspond to points on the red curves that are close to the minimum of these curves.

\begin{figure}[ht]
 \centering
 {\includegraphics[scale=0.32]{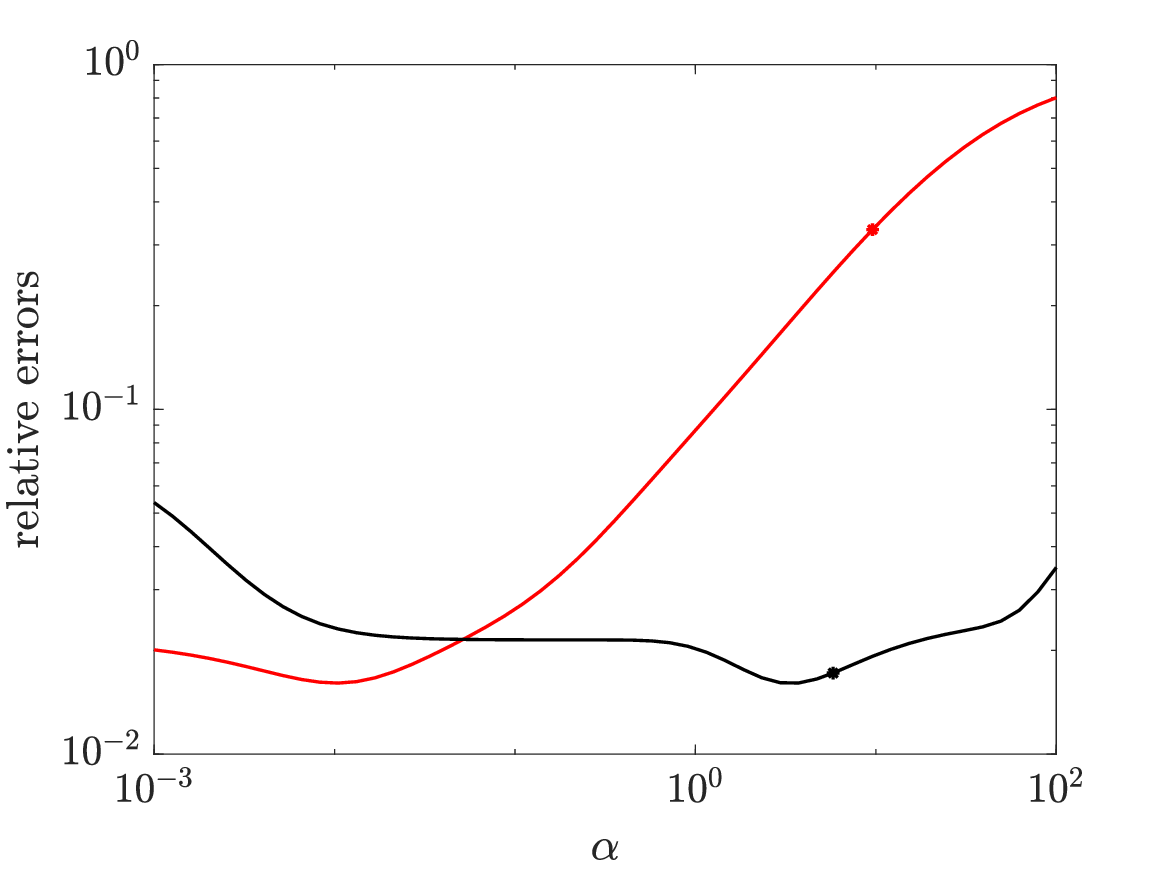}}
 \hspace{0.5cm}
 {\includegraphics[scale=0.32]{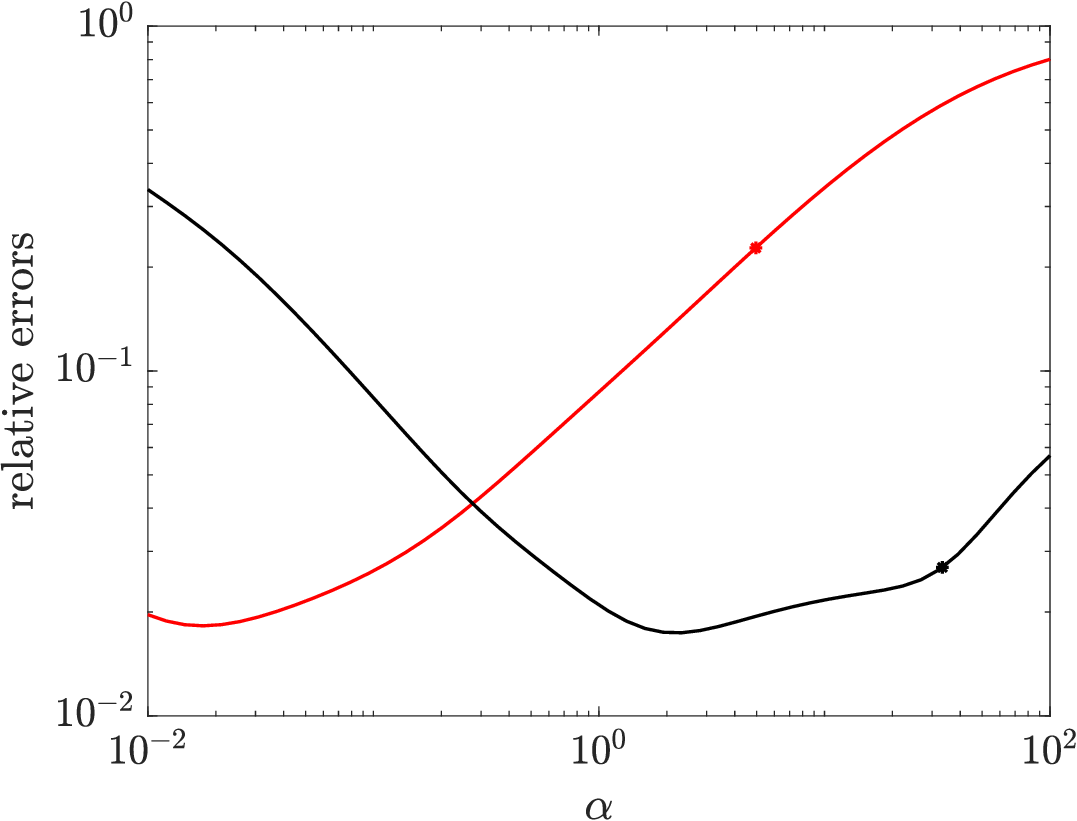}}
 \caption{Example~\ref{ex1} - Relative error in approximate solutions computed 
 by AT (red) and iAT (black) when varying $\alpha$ for $n=1000$ and $\xi=0.01$, i.e. $\delta = 1\%$. The points marked by $\ast$ correspond to the value of $\alpha$ 
 determined by solving~\eqref{alpha_iAT} for iAT and as in~\cite{ramlau2019} for AT. (Left) $\ell=10$, $i=200$, (Right) $\ell=30$, 
 $i=100$.}\label{fig2}
\end{figure}

\begin{figure}[ht]
\centering
{\includegraphics[scale=0.32]{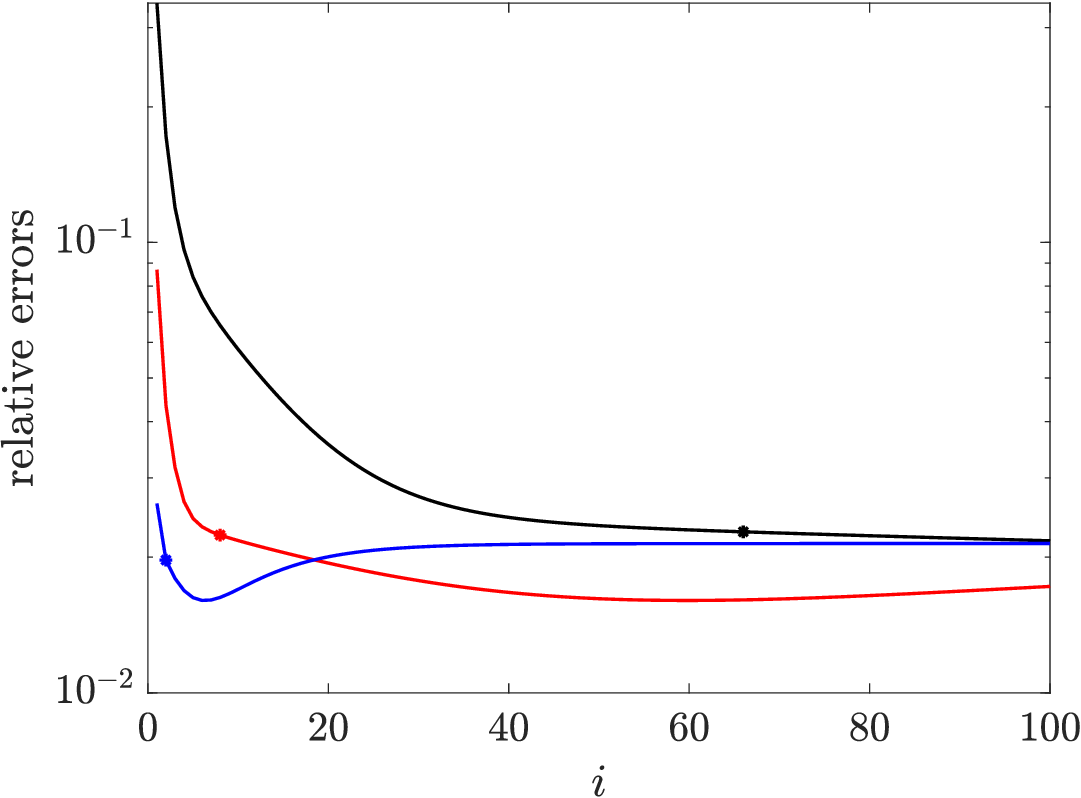}}
\hspace{0.5cm}
{\includegraphics[scale=0.32]{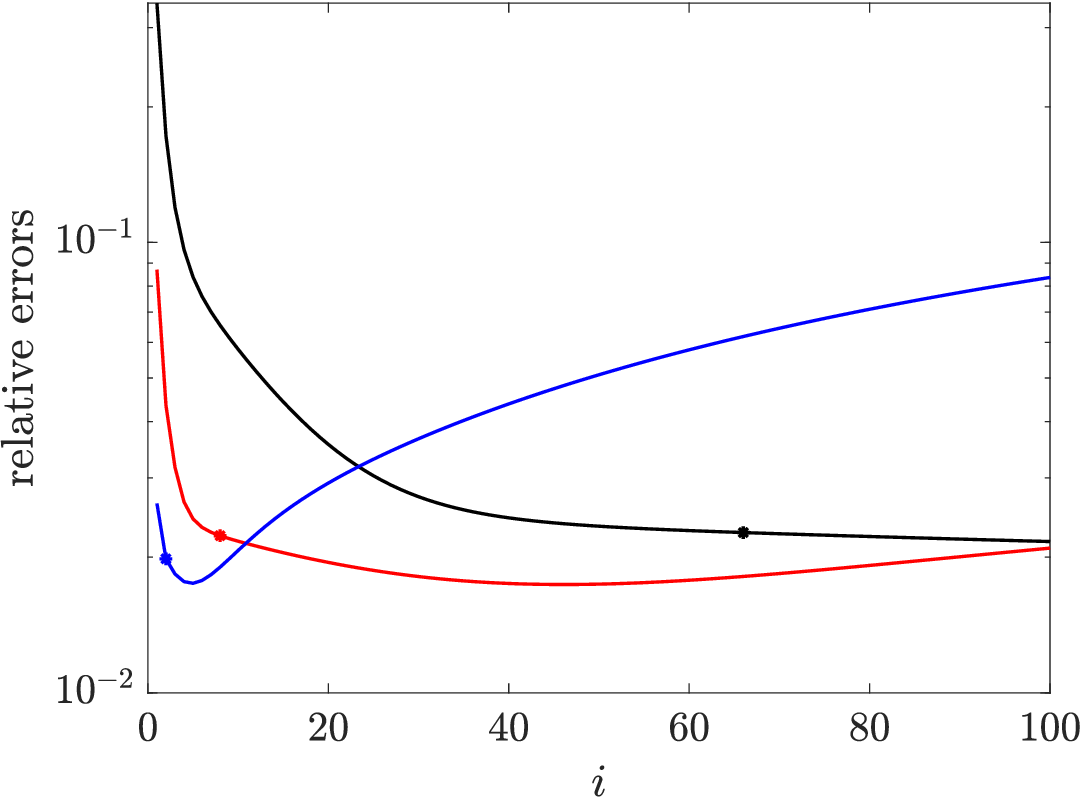}}
\caption{Example~\ref{ex1} - Relative error in approximate solutions computed by iAT as a 
function of $i$ for $\alpha=10$ (black), $\alpha=1$ (red), and $\alpha=0.1$ (blue), for
$n=1000$ and $\xi=0.01$, i.e. $\delta = 1\%$. The points marked by $*$ show when the discrepancy principle is 
satisfied; see Table~\ref{tab2}. (Left) $\ell=10$, (Right) $\ell=30$.}\label{fig3}
\end{figure}

We turn to the behaviour of Algorithm~\ref{algo:itarntik} for fixed values of $\alpha$ and
determine the number of iterations $i$ by the discrepancy principle, i.e., we let $i$ be 
the smallest index such that 
\begin{equation*}
     \|T_nx_{\alpha,n,i}^{\delta,\ell}-y_n^{\delta}\|_2 \leq \delta.
\end{equation*}
Figure~\ref{fig3} shows the relative error in approximate solutions computed with
Algorithm~\ref{algo:itarntik} as a function of $i$ for three values of $\alpha$ and 
$\ell=10$ (Figure~\ref{fig3} (a)) and $\ell=30$ (Figure~\ref{fig3} (b)). When $\alpha$ is 
small, it suffices to carry out a few iterations to satisfy the discrepancy principle at the
cost of minor instability when increasing $i$. This behavior does not change much when 
increasing $\ell$. This is illustrated in Figure~\ref{fig3} (b) for $\ell=30$. The figure 
shows semi-convergence, i.e., that the error first decreases and subsequently increases
when $i$ grows, to be more pronounced when $\alpha>0$ is small. Note that the discrepancy 
principle terminates the iterations for the same $i$-values for both $\ell$-values. Some 
relative errors in the computed approximate solutions are displayed in Table~\ref{tab2}. 
We note that when $i=1$, the iterated Tikhonov method simplifies to a noniterated 
Tikhonov method. The latter differs from the AT method in~\cite{ramlau2019} in that 
$\alpha$ is chosen differently. The relative errors in Table \ref{tab2} are close to the 
smallest ones of Table~\ref{tab1}, but are achieved with fewer $i$-iterations.

\begin{table}
\caption{Example~\ref{ex1} - Relative error in approximate solutions computed
by iAT for fixed values of $\alpha$, and $i$ determined by the discrepancy 
principle for $n=1000$ and $\xi=0.01$, i.e. $\delta = 1\%$.}\label{tab2}
  \begin{tabular}{cccccc}
 \toprule%
 & \multicolumn{3}{c}{$\ell=10$} & & $\ell=30$ \\
 \cmidrule{2-3}
 \cmidrule{5-6}
  $\alpha$ & $i$ & $\|x_n^{\dagger}-x_{\alpha,n,i}^{\delta,\ell}\|_2/\|x_n^{\dagger}\|_2$ & & $i$ & $\|x_n^{\dagger}-x_{\alpha,n,i}^{\delta,\ell}\|_2/\|x_n^{\dagger}\|_2$\\
  \midrule
  10 & 66 & $2.28\cdot 10^{-2}$ & & 66 & $2.27\cdot 10^{-2}$\\
  5 & 34 & $2.27\cdot 10^{-2}$ & & 34 & $2.26\cdot 10^{-2}$\\
  1 & 8 & $2.24\cdot 10^{-2}$ & & 8 & $2.23\cdot 10^{-2}$\\
  0.5 & 5 & $2.18\cdot 10^{-2}$ & & 5 & $2.18\cdot 10^{-2}$\\
  0.1 & 2 & $1.97\cdot 10^{-2}$ & & 2 & $1.98\cdot 10^{-2}$\\
  0.01 & 1 & $1.61\cdot 10^{-2}$ & & 1 & $1.96\cdot 10^{-2}$\\
  \botrule
  \end{tabular}
\end{table}

\begin{figure}[ht]
\centering
\includegraphics[scale=0.35]{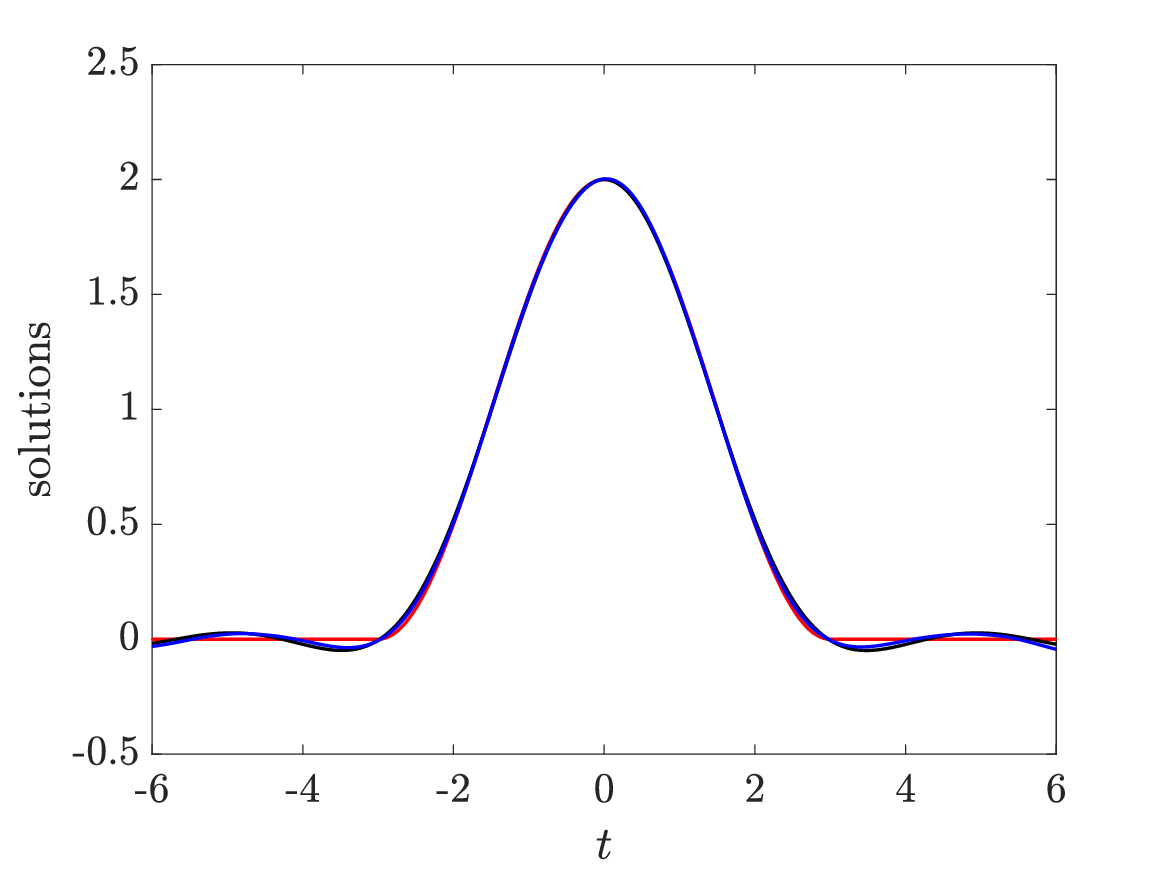}
\caption{Example~\ref{ex1} - Exact solution $x_n^{\dagger}$ (red) and approximate solutions $x_{\alpha,n,i}^{\delta,\ell}$ (black) computed by iAT with $i=100$ and $\alpha$ determined by solving~\eqref{alpha_iAT}, and approximate 
solution $x_{0.1,n,2}^{\delta,\ell}$ (blue) computed by iAT using the discrepancy principle for $n=1000$, $\ell=30$,
and $\xi=0.01$, i.e. $\delta = 1\%$.}\label{fig4}
\end{figure}

Figure~\ref{fig4} depicts the exact solution $x_n^{\dagger}$ as well as approximate 
solutions determined by Algorithm~\ref{algo:itarntik} for $\ell=30$. The approximate 
solution plotted in black is computed with Algorithm~\ref{algo:itarntik} (iAT) with 
$i=100$, and the approximate solution plotted in blue is determined by setting 
$\alpha=0.1$ and terminating the $i$-iterations with the discrepancy principle. Note that 
both computed approximate solutions are accurate approximations of $x_n^{\dagger}$.
\end{exmp}

\begin{exmp}\label{ex2}
We consider the Fredholm integral equation of the first kind discussed by~Baart~\cite{baart1982},
\begin{equation*}
  \int_{-6}^{6}\kappa(s,t)x(t)dt=y(s),\qquad 0\leq s\leq\frac{\pi}{2}, 
\end{equation*}
where
\begin{equation*}
   x(t)=\sin(t),\quad\kappa(s,t)=\exp(s\cos(t)),\quad\text{and}\quad y(s)=2\sinh(s)/s.
\end{equation*}
The discretization is determined with the MATLAB function \texttt{baart} from 
\cite{hansen2007}, which gives a nonsymmetric matrix $T_n\in\R^{n\times n}$ and the true 
solution $x_n^{\dagger}\in\R^n$.

\begin{table}
\caption{Example~\ref{ex2} - Relative error in approximate solutions computed by iAT and $\alpha$ determined by solving $\eqref{alpha_iAT}$ for different
values of $\ell$, with $n=1000$ and $\xi=0.01$, i.e. $\delta = 1\%$. The AT method is applied 
with the parameter $\alpha$ determined as in~\cite{ramlau2019}.}\label{tab4}
 \begin{tabular}{cccccc}
 \toprule%
 & \multicolumn{3}{c}{iAT} & & AT \\
 \cmidrule{2-4}
 \cmidrule{6-6}
 $\ell$ & $i$ & $\alpha$ & $\|x_n^{\dagger}-x_{\alpha,n,i}^{\delta,\ell}\|_2/
 \|x_n^{\dagger}\|_2$ & & $\|x_n^{\dagger}-x_{\alpha,n}^{\delta,\ell}\|_2/
 \|x_n^{\dagger}\|_2\ (\alpha)$\\
 \midrule
 \multirow{3}*{3} & 1 & $1.23\cdot 10^{0\phantom{-}}$ & $5.87\cdot 10^{-1}$ & & 
 \multirow{3}*{$7.48\cdot 10^{-1} \ (8.28)$}\\
 & 200 & $5.86\cdot 10^{0\phantom{-}}$ & $3.67\cdot 10^{-1}$ & \\
 & 500 & $2.03\cdot 10^{0\phantom{-}}$ & $2.74\cdot 10^{-1}$ &\\
 \midrule
 \multirow{3}*{6} & 1 & $1.27\cdot 10^{-1}$ & $3.67\cdot 10^{-1}$ & &
 \multirow{3}*{$4.25\cdot 10^{-1}\ (2.78\cdot 10^{-1})$}\\
  & 200 & $5.86\cdot 10^{0\phantom{-}}$ & $3.04\cdot 10^{-1}$ & \\
  & 500 & $2.03\cdot 10^{0\phantom{-}}$ & $1.84\cdot 10^{-1}$ &\\
 \midrule
 \multirow{3}*{9} & 1 & $5.44\cdot 10^{-2}$ & $3.32\cdot 10^{-1}$ & & 
 \multirow{3}*{$3.37\cdot 10^{-1}\ (6.37\cdot 10^{-2})$}\\
  & 200 & $5.86\cdot 10^{0\phantom{-}}$ & $3.05\cdot 10^{-1}$ & \\
  & 500 & $2.03\cdot 10^{0\phantom{-}}$ & $1.90\cdot 10^{-1}$ &\\
  \botrule
 \end{tabular}
\end{table}

\begin{table}
\caption{Example~\ref{ex2} - Relative error in approximate solutions computed by iAT and $\alpha$ determined by solving $\eqref{alpha_iAT}$ for different
values of $\ell$, with $n=1000$ and $\xi=0.001$, i.e. $\delta = 0.1\%$. The AT method is applied 
with the parameter $\alpha$ determined as in~\cite{ramlau2019}.}\label{tab5}
  \begin{tabular}{cccccc}
 \toprule%
 & \multicolumn{3}{c}{iAT} & & AT \\
 \cmidrule{2-4}
 \cmidrule{6-6}
  $\ell$ & $i$ & $\alpha$ & $\|x_n^{\dagger}-x_{\alpha,n,i}^{\delta,\ell}\|_2/
  \|x_n^{\dagger}\|_2$ & & $\|x_n^{\dagger}-x_{\alpha,n}^{\delta,\ell}\|_2/
  \|x_n^{\dagger}\|_2\ (\alpha)$\\
  \midrule
  \multirow{3}*{3} & 1 & $6.51\cdot 10^{-1}$ & $5.17\cdot 10^{-1}$ & & 
  \multirow{3}*{$7.14\cdot 10^{-1}\ (5.48)$}\\
  & 500 & $2.03\cdot 10^{0\phantom{-}}$ & $1.28\cdot 10^{-1}$ &\\
  & 1000 & $1.43\cdot 10^{0\phantom{-}}$ & $4.60\cdot 10^{-2}$ &\\
  \midrule
  \multirow{3}*{6} & 1 & $7.49\cdot 10^{-2}$ & $3.42\cdot 10^{-1}$ & &
  \multirow{3}*{$4.09\cdot 10^{-1}\ (2.32\cdot 10^{-1})$}\\
  & 500 & $2.03\cdot 10^{0\phantom{-}}$ & $1.80\cdot 10^{-1}$ &\\
  & 1000 & $1.43\cdot 10^{0\phantom{-}}$ & $1.56\cdot 10^{-1}$ &\\
  \midrule
  \multirow{3}*{9} & 1 & $2.29\cdot 10^{-3}$ & $1.90\cdot 10^{-1}$ & & 
  \multirow{3}*{$2.15\cdot 10^{-1}\ (4.23\cdot 10^{-3})$}\\
  & 500 & $1.45\cdot 10^{0\phantom{-}}$ & $1.74\cdot 10^{-1}$ &\\
  & 1000 & $1.43\cdot 10^{0\phantom{-}}$ & $1.65\cdot 10^{-1}$ &\\
  \botrule
  \end{tabular}
\end{table}

Table~\ref{tab4} shows the relative error in approximate solutions computed by the AT and
iAT methods for $n=1000$ and noise level $\delta=1\%$. The singular values of the matrix
$T_n$, when ordered in nonincreasing order, decreases very rapidly with increasing index. 
Therefore, it is not meaningful to choose $\ell$ larger than $9$. Algorithm 
\ref{algo:itarntik} (iAT) can be seen to yield approximate solutions with smaller relative 
error than Algorithm \ref{algo:AT} (AT), in particular for small values of $\ell$. Table
\ref{tab5} differs from Table \ref{tab4} only in that the noise level is $\delta=0.1\%$.

\begin{figure}[ht]
\centering
{\includegraphics[scale=0.31]{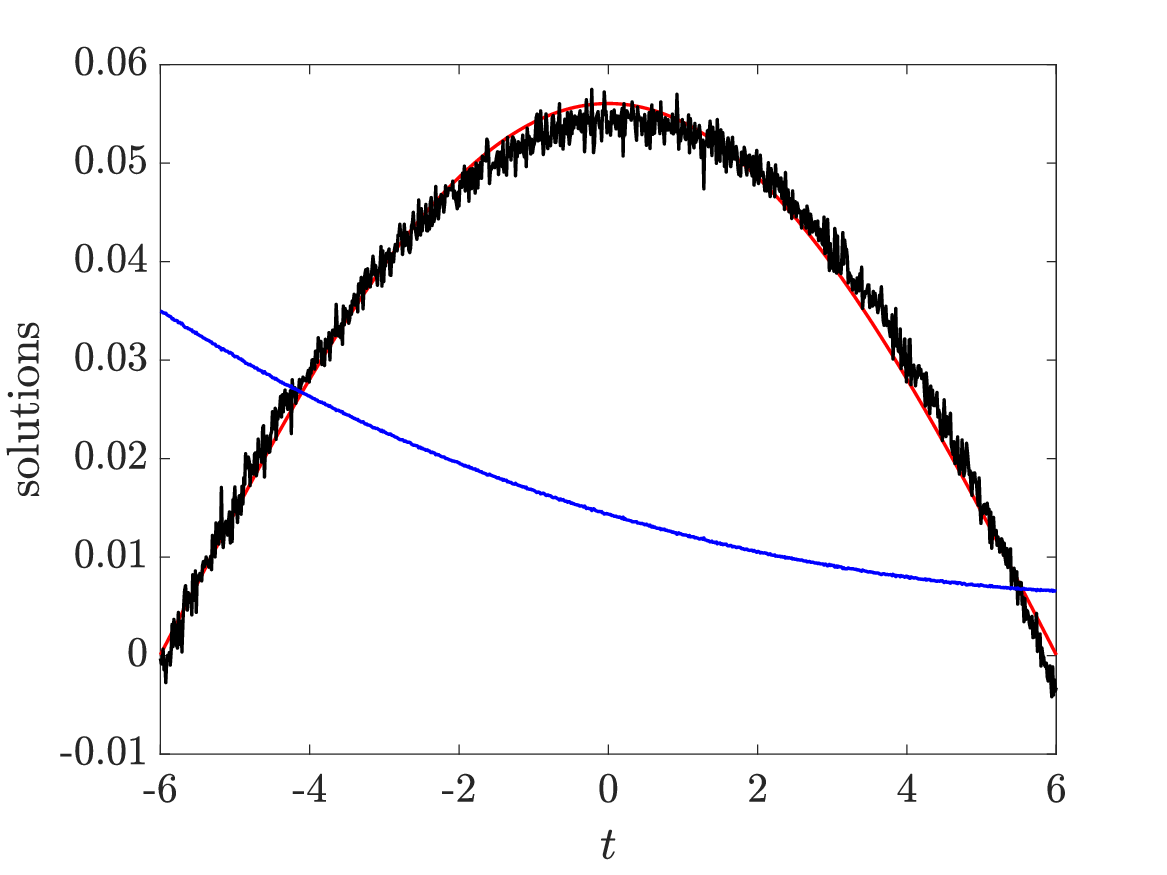}}
\hspace{0.5cm}
{\includegraphics[scale=0.31]{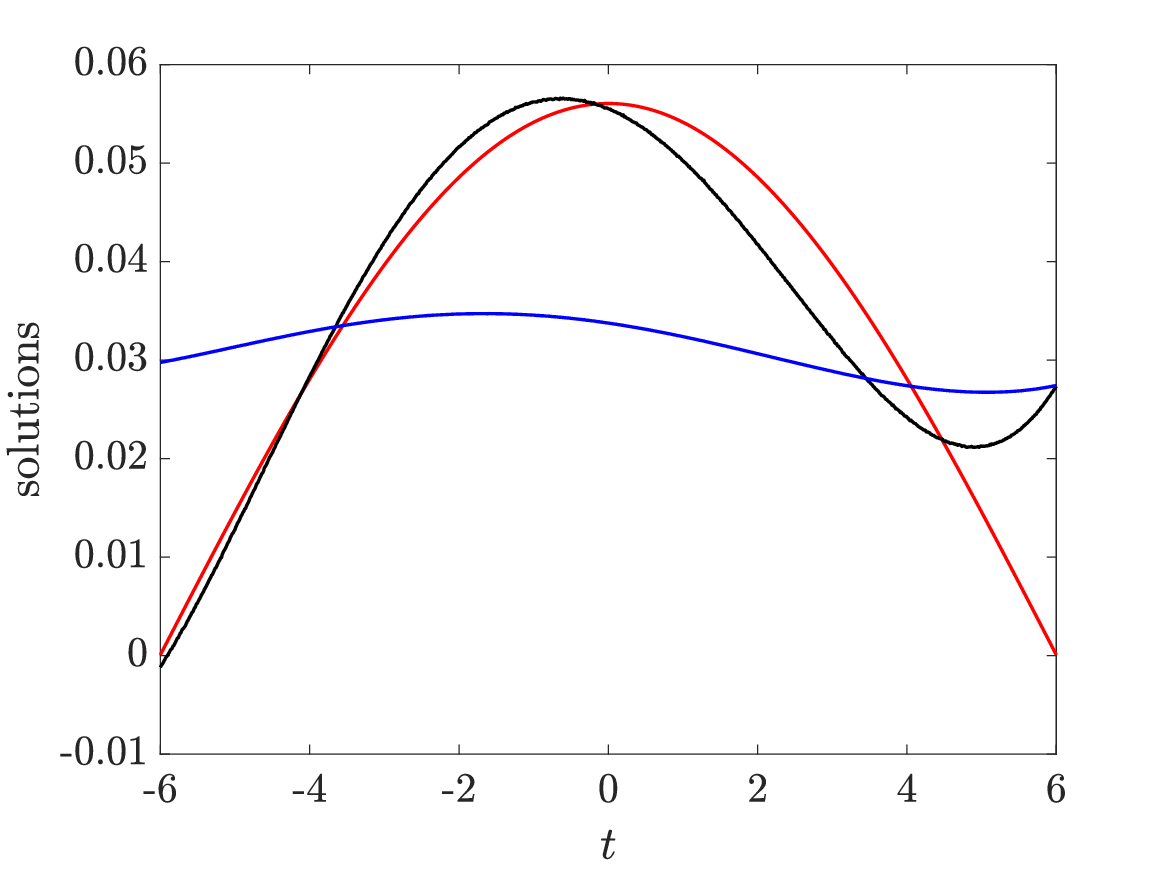}}
\caption{Example~\ref{ex2} - Exact solution $x^{\dagger}_n$ (red) and approximate 
solutions $x^{\delta,\ell}_{\alpha,n,i}$ computed by iAT (black) with $i=1000$ and $\alpha$ determined by solving~\eqref{alpha_iAT} and 
$x^{\delta,\ell}_{\alpha,n}$ computed by AT (blue) with $\alpha$ determined as in~\cite{ramlau2019}, for $n=1000$, $\xi=0.001$, i.e. $\delta = 0.1\%$. 
(Left) $\ell=3$, (Right) $\ell=6$.}\label{fig6}
\end{figure}

\begin{figure}[ht]
\centering
{\includegraphics[scale=0.32]{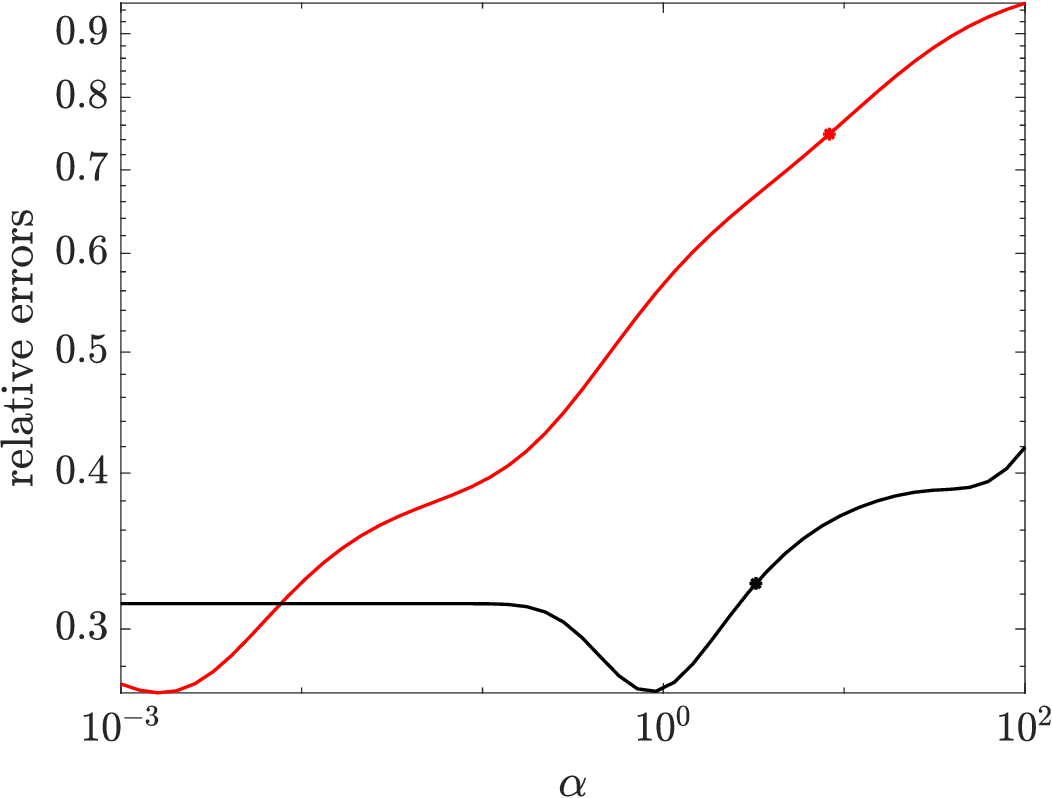}}
\hspace{0.5cm}
{\includegraphics[scale=0.32]{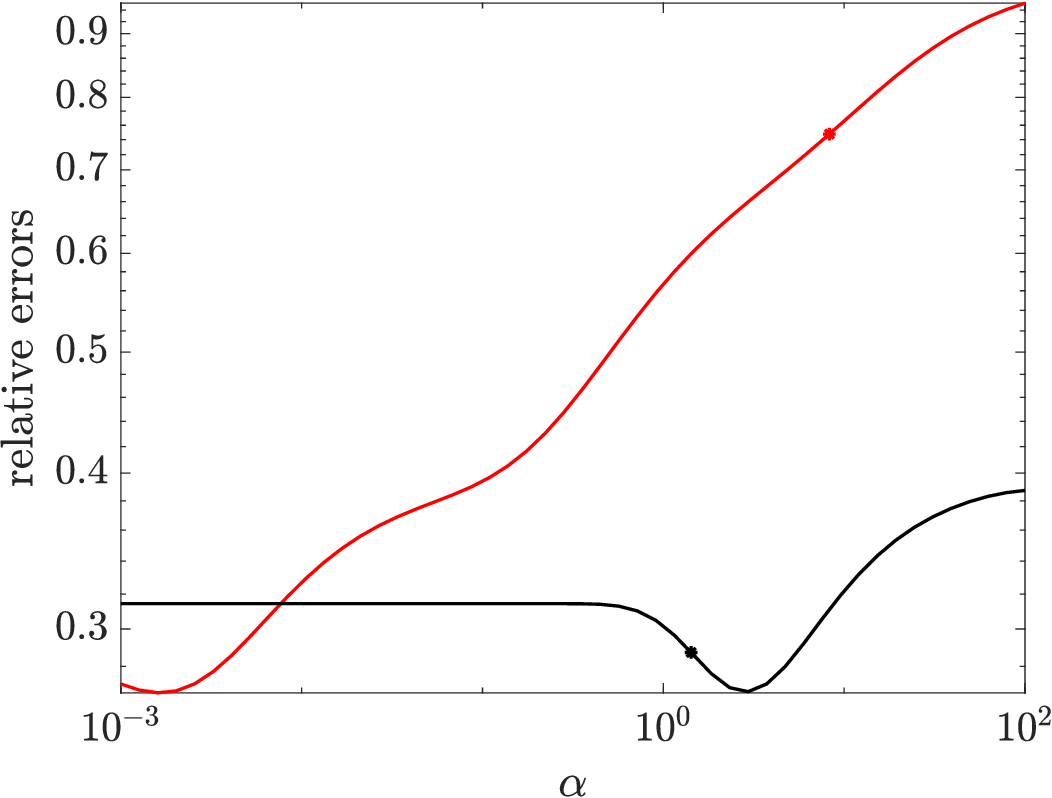}}
\caption{Example~\ref{ex2} - Relative error in approximate solutions computed 
by AT (red) and iAT (black) when varying $\alpha$ for 
$n=1000$, $\ell=3$ and $\xi=0.01$, i.e. $\delta = 1\%$. The points marked by $\ast$ correspond to the value of $\alpha$ 
determined by solving~\eqref{alpha_iAT} for iAT and as in~\cite{ramlau2019} for AT. (Left) $i=300$, (Right) 
$i=1000$.}
\label{fig7}
\end{figure}

Figure~\ref{fig6} depicts the exact solution and computed approximate solutions for two 
values of $\ell$ and noise level $0.1\%$. These plots illustrate the improved quality of 
the computed solutions determined by the iAT method when compared with approximate 
solutions determined by the AT method. Figure~\ref{fig7} displays the behavior of the 
relative error in the computed approximate solutions when varying the parameter $\alpha$.
The value of $\alpha$ for the iAT method, which is determined by solving~\eqref{alpha_iAT}, corresponds to a point that is closer to the 
minimum of the relative error than the point associated with the parameter value for the
AT with $\alpha$ determined as in~\cite{ramlau2019}. 

\begin{figure}[ht]
\centering
{\includegraphics[scale=0.32]{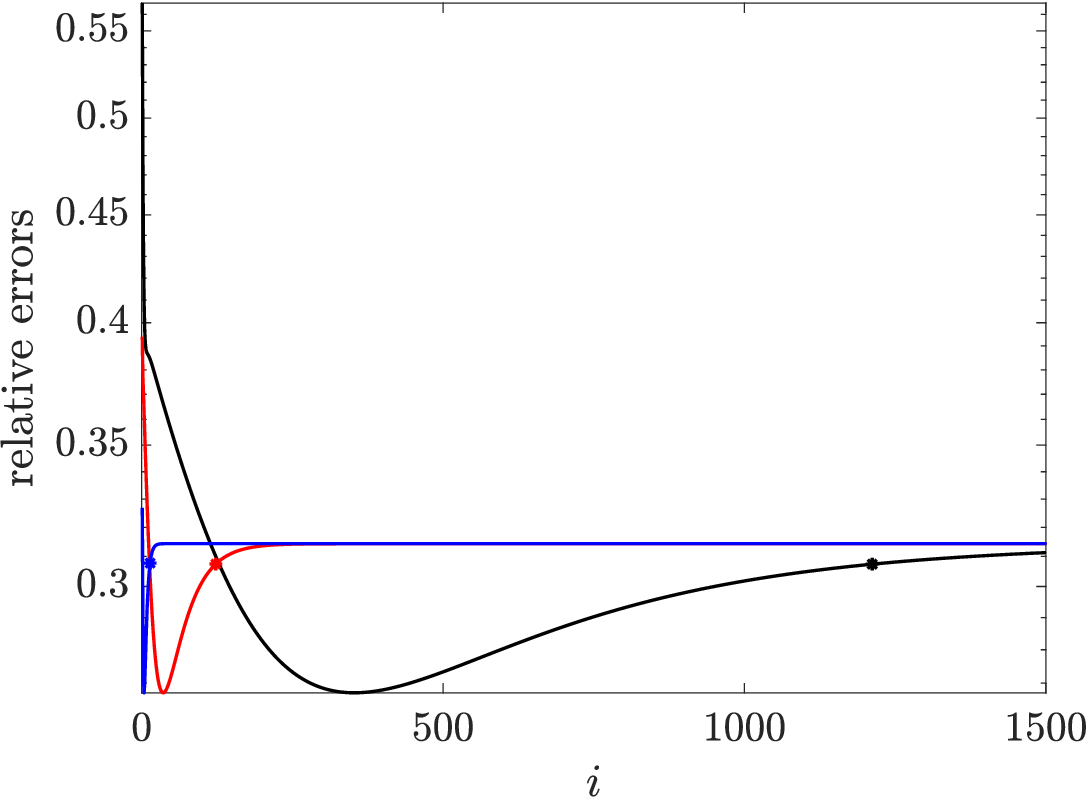}}
\hspace{0.5cm}
{\includegraphics[scale=0.32]{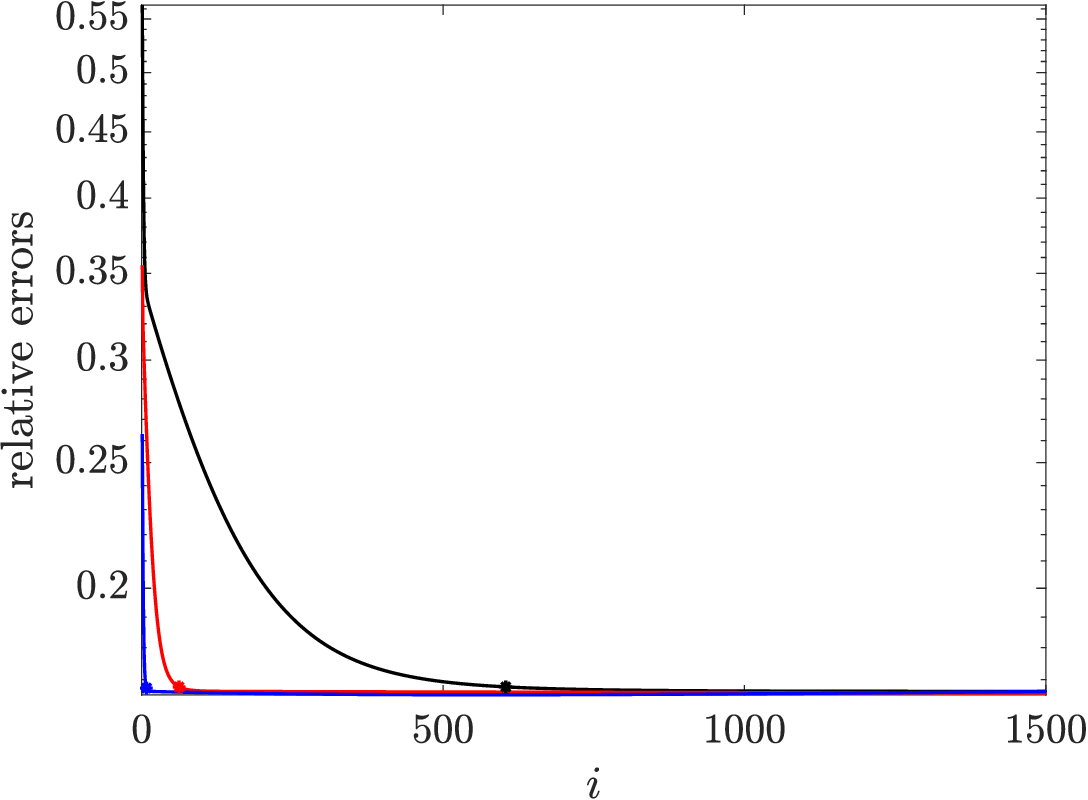}}
\caption{Example~\ref{ex2} - Relative error in approximate solutions computed by iAT as a 
function of $i$ for $\alpha=1$ (black), $\alpha=0.1$ (red), and $\alpha=0.01$ (blue), for
$n=1000$ and $\xi=0.01$, i.e. $\delta = 1\%$. The points marked by $*$ show when the discrepancy principle is 
satisfied; see Table~\ref{tab6}. (Left) $\ell=3$, 
(Right) $\ell=9$.}\label{fig8}
\end{figure}

We turn to the relative error of computed approximate solutions as a function of the 
iteration number $i$ for fixed values of $\alpha$. Figure~\ref{fig8} shows the relative 
error for three values of $\alpha$ and varying $i$ for two values of $\ell$. Note that, as
expected, a larger $\ell$ gives a smaller error. Points marked by $*$ on the curves 
show the smallest $i$-values for which the discrepancy principle is satisfied. We note 
that the discrepancy principle is satisfied for a small $i$-value when $\alpha$ is small. 

\begin{table}
\caption{Example~\ref{ex2} - Relative error in approximate solutions computed
by iAT for fixed values of $\alpha$, and $i$ determined by the discrepancy 
principle for $n=1000$ and $\xi=0.01$, i.e. $\delta = 1\%$.}\label{tab6}
  \begin{tabular}{cccccc}
 \toprule%
 & \multicolumn{3}{c}{$\ell=3$} & & $\ell=9$ \\
 \cmidrule{2-3}
 \cmidrule{5-6}
  $\alpha$ & $i$ & $\|x_n^{\dagger}-x_{\alpha,n,i}^{\delta,\ell}\|_2/
  \|x_n^{\dagger}\|_2$ & & $i$ & $\|x_n^{\dagger}-x_{\alpha,n,i}^{\delta,\ell}\|_2/
  \|x_n^{\dagger}\|_2$\\
  \midrule
  1 & 1212 & $3.07\cdot 10^{-1}$ & & 604 & $1.68\cdot 10^{-1}$\\
  0.1 & 123 & $3.07\cdot 10^{-1}$ & & 62 & $1.68\cdot 10^{-1}$\\
  0.01 & 14 & $3.08\cdot 10^{-1}$ & & 8 & $1.67\cdot 10^{-1}$\\
  \botrule
  \end{tabular}
\end{table}

\begin{figure}[ht]
\centering
\includegraphics[scale=0.35]{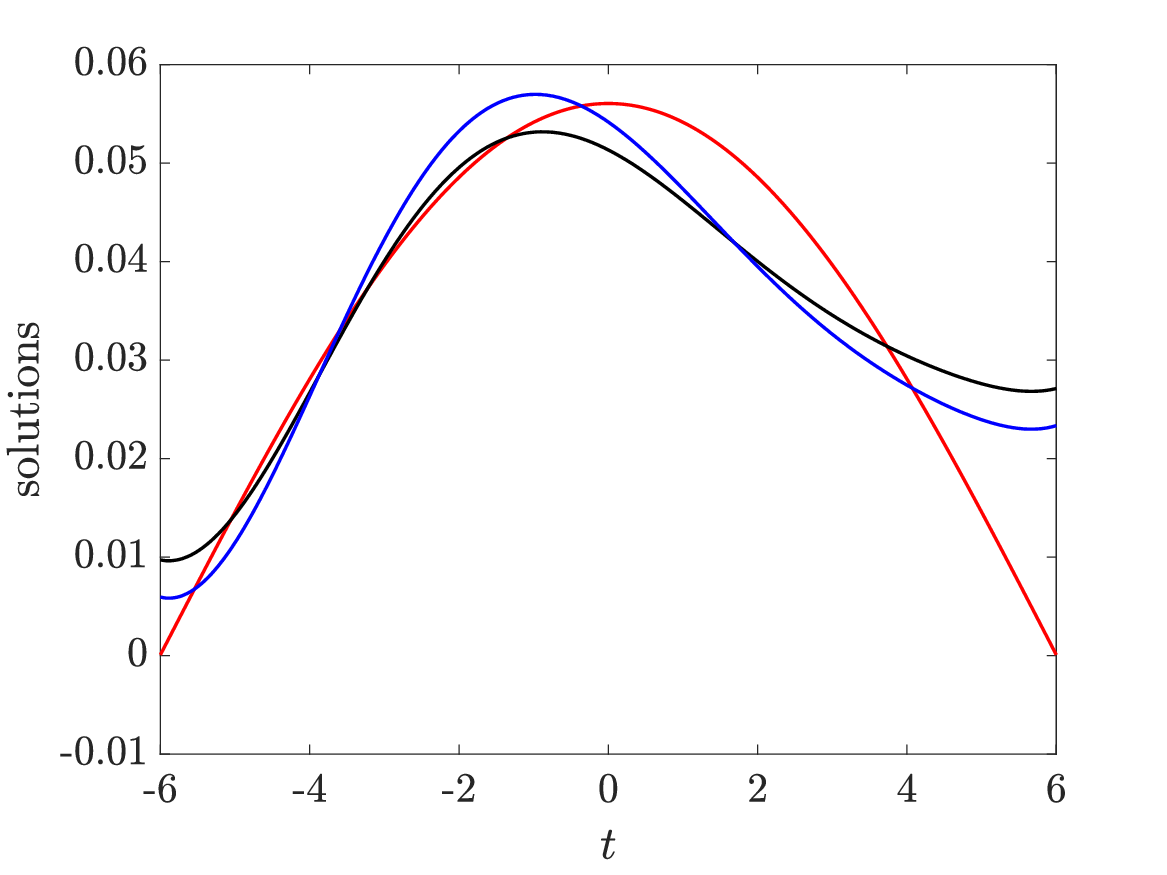}
\caption{Example~\ref{ex2} - Exact solution $x_n^{\dagger}$ (red) and approximate solutions $x_{\alpha,n,i}^{\delta,\ell}$ (black) computed by iAT with $i=500$ and $\alpha$ determined by solving~\eqref{alpha_iAT}, and approximate 
solution $x_{0.1,n,62}^{\delta,\ell}$ (blue) computed by iAT using the discrepancy principle for $n=1000$, $\ell=9$,
and $\xi=0.01$, i.e. $\delta = 1\%$.}\label{fig9}
\end{figure}

Table~\ref{tab6} displays the relative error when $i$ is determined by the discrepancy 
principle for three values of $\alpha$. The relative error in the computed approximate
solution is essentially independent of the value of $\alpha$. For $\ell=9$ the results are
comparable to those of Table~\ref{tab5}. Figure~\ref{fig9} displays the exact solution for
$n=1000$ (red curve) and computed solutions for $\ell=9$. The black curve displays the 
computed solution for $i=500$ for $\alpha$ determined by solving~\eqref{alpha_iAT}; the blue 
curve depicts the computed solution obtained by fixing $\alpha=0.1$ and determining $i$ 
with the discrepancy principle.
\end{exmp}

\begin{exmp}\label{ex3}
This example is concerned with a digital image deblurring problem. We use the function 
\texttt{blur} from~\cite{hansen2007} with default parameters to determine an $n^2\times n^2$ symmetric block 
Toeplitz matrix with Toeplitz blocks, $T_{n^2}$, that models blurring of an image that is
represented by $n\times n$ pixels. The blur is determined by a space-invariant Gaussian 
point spread function. The true image is represented by the vector 
$x_{n^2}^{\dagger}\in\R^{n^2}$ for $n=30$. This image is shown in Figure~\ref{fig2D} (a).

\begin{table}
\caption{Example~\ref{ex3} - Relative error in approximate solutions computed by iAT and $\alpha$ determined by solving $\eqref{alpha_iAT}$ for $\ell=300$, $n^2=900$ and $\xi=0.01$, i.e. $\delta = 1\%$.}\label{tab7}
  \begin{tabular}{ccc}
  \toprule
  $i$ & $\alpha$ & $\|x_{n^2}^{\dagger}-x_{\alpha,n^2,i}^{\delta,\ell}\|_2/\|x_{n^2}^{\dagger}\|_2$\\
  \midrule
  1 & $2.44\cdot 10^1$ & $9.72\cdot 10^{-1}$\\
  100 & $3.36\cdot 10^1$ & $3.68\cdot 10^{-1}$\\
  200 & $5.82\cdot 10^0$ & $1.37\cdot 10^{-1}$\\
  300 & $3.24\cdot 10^0$ & $9.01\cdot 10^{-2}$\\
  400 & $2.42\cdot 10^0$ & $7.56\cdot 10^{-2}$\\
  500 & $2.02\cdot 10^0$ & $6.83\cdot 10^{-2}$\\
  1000 & $1.42\cdot 10^0$ & $6.05\cdot 10^{-2}$\\
  \botrule
  \end{tabular}
\end{table}

Table~\ref{tab7} displays relative restoration errors for the computed approximate 
solutions with Algorithm \ref{algo:itarntik} (iAT) for the noise level $\delta = 1\%$. To 
satisfy inequality~\eqref{condEC2}, we need to let $\ell>268$; any larger value of $\ell$
gives very similar results, in particular for a large number of iterations $i$. Therefore,
we set $\ell=300$.

Figure~\ref{fig2D} (b) shows the approximate solution computed with iAT for $\ell=300$ and
$i=500$; $\alpha$ is determined by solving equation~\eqref{alpha_iAT}. For the AT method, we 
need at least $\ell=756$ to satisfy~\eqref{condEC2} with the choices $C=1$ and $E=3\|x_{n^2}^{\dagger}\|_2$. Figure~\ref{fig2D} (c) shows the 
approximate solution computed by AT for $\ell=800$ and for $\alpha$ determined as in 
\cite{ramlau2019}.

\begin{figure}[ht]
\centering
{\includegraphics[scale=0.3]{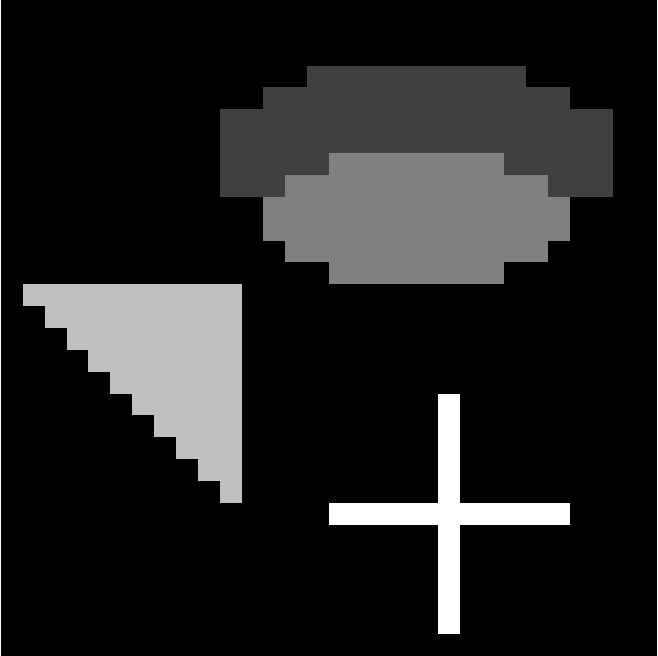}} 
\hspace{1cm}
{\includegraphics[scale=0.365]{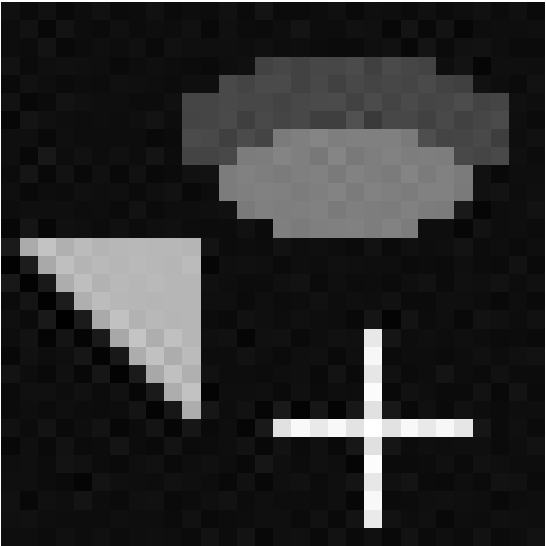}}
\hspace{1cm}
{\includegraphics[scale=0.365]{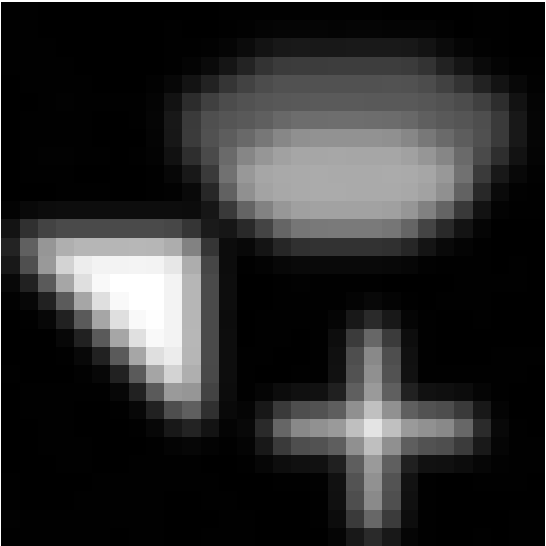}}
\caption{Example~\ref{ex3} - Exact solution $x^{\dagger}_n$ (Left) and approximate 
solutions $x^{\delta,\ell}_{\alpha,n,i}$ computed by iAT (Center) with $\ell=300$, $i=500$ and $\alpha$ determined by solving~\eqref{alpha_iAT} and 
$x^{\delta,\ell}_{\alpha,n}$ computed by AT (Right) with $\ell=800$ and $\alpha$ determined as in~\cite{ramlau2019}, for $n^2=900$, $\xi=0.01$, i.e. $\delta = 1\%$.}\label{fig2D}
\end{figure}

\begin{figure}[ht]
\centering
{\includegraphics[scale=0.32]{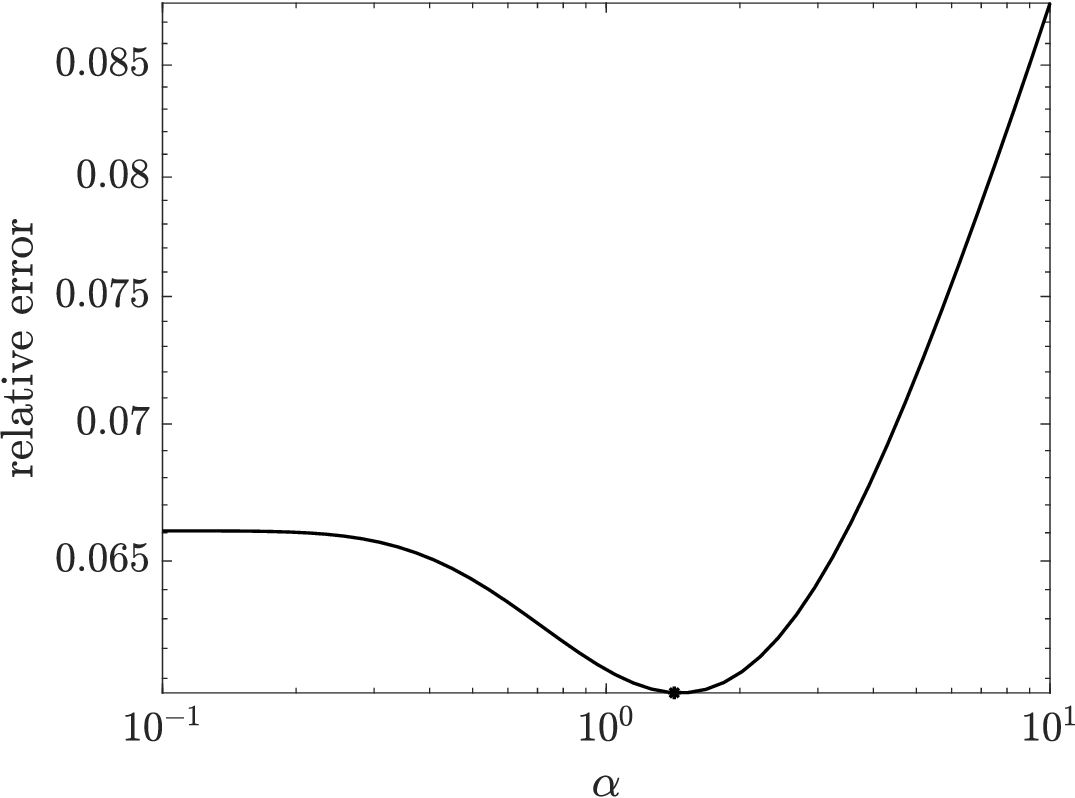}} 
\hspace{0.5cm}
{\includegraphics[scale=0.32]{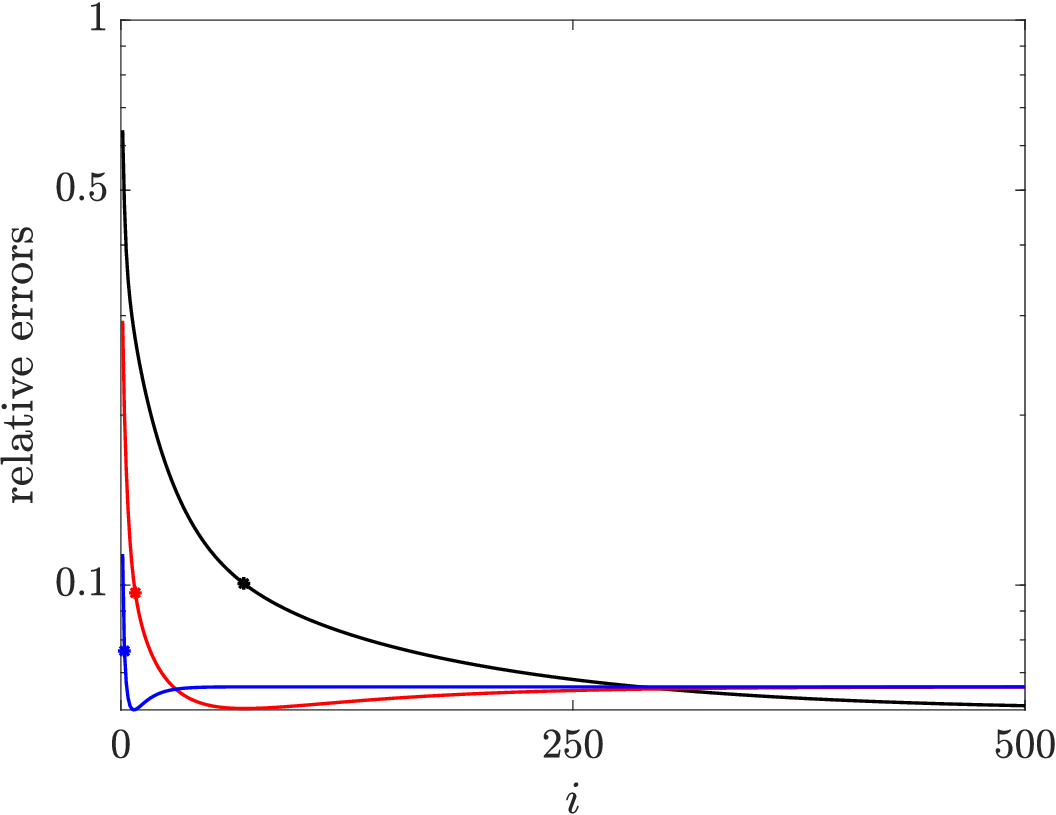}}
\caption{Example~\ref{ex3} - Relative error in approximate solutions computed by iAT for $n^2=900$, $\ell=300$ and $\xi=0.01$, i.e. $\delta = 1\%$. (Left) Relative error varying $\alpha$, where the point marked by $\ast$ correspond to the value of $\alpha$ 
determined by solving~\eqref{alpha_iAT} with $i=1000$, (Right) Relative error as a function of $i$ for $\alpha=1$ (black), $\alpha=0.1$ (red), and $\alpha=0.01$ (blue). The points marked by $*$ show when the discrepancy principle is 
satisfied; see Table~\ref{tab8}.}\label{fig11}
\end{figure}

\begin{table}
\centering
  \begin{tabular}{ccc}
  \toprule
  $\alpha$ & $i$ & $\|x_{n^2}^{\dagger}-x_{\alpha,n^2,i}^{\delta,\ell}\|_2/\|x_{n^2}^{\dagger}\|_2$\\
  \midrule
  5 & 333 & $1.01\cdot 10^{-1}$\\
  1 & 68 & $1.01\cdot 10^{-1}$\\
  0.5 & 35 & $9.99\cdot 10^{-2}$\\
  0.1 & 8 & $9.69\cdot 10^{-2}$\\
  0.05 & 5 & $9.08\cdot 10^{-2}$\\
  0.01 & 2 & $7.64\cdot 10^{-2}$\\
  0.001 & 2 & $6.43\cdot 10^{-2}$\\
  \botrule
  \end{tabular}
\caption{Example~\ref{ex3} - Relative error in approximate solutions computed
by iAT for fixed values of $\alpha$, and $i$ determined by the discrepancy 
principle for $n^2=900$ and $\xi=0.01$, i.e. $\delta = 1\%$.}\label{tab8}
\end{table}

Figure~\ref{fig11} displays the relative error of the computed solutions when varying the 
parameters $\alpha$ and $i$. Specifically, Figure~\ref{fig11} (a) is obtained by varying 
$\alpha$ for $i=500$. The $*$ on the graph corresponds to the value of the parameter 
$\alpha$ determined by solving equation~\eqref{alpha_iAT}. Similarly, Figure~\ref{fig11} (b) is 
obtained by varying $i$ for three fixed values of $\alpha$. The $*$s mark $i$-values 
for which the discrepancy principle is satisfied. The robustness of the discrepancy 
principle is illustrated by Table~\ref{tab8}, where for each $\alpha$, the number of 
iterations $i$ is determined by the discrepancy principle. The table displays the relative
error in the computed approximate solutions.

We remark that the system matrix $T_{n^2}$ in this example is symmetric. Therefore, the 
symmetric Lanczos algorithm could be applied instead of the Arnoldi algorithm. However,
computed examples reported in~\cite{AR} show the latter algorithm to yield higher
accuracy because the computed Lanczos vectors generated by the Lanczos algorithm may be 
far from orthogonal due to propagated and amplified round-off errors.
\end{exmp}

\section{Conclusion and extensions}\label{sec:end}
The paper presents a convergence analysis for iterated Tikhonov regularization based on
the Arnoldi decomposition. A new approach to choosing the regularization parameter is 
proposed. Computed results show iterated Tikhonov regularization with the proposed
parameter choice to yield more accurate approximate solutions than the non-iterative
method studied in~\cite{ramlau2019}. It would be interesting to extend this work to 
Tikhonov regularization in general form. Applications of such problems are described in
e.g.,~\cite{bianchi2022,bianchi2015,bianchi2023graph,bianchi2025data,bianchi2023uniformly,
gazzola2015krylov,huang19,bianchi2017generalized}. It also may be interesting to compare the choice of 
regularization parameter of this paper with other parameter choices, such as those 
discussed in~\cite{kindermann2011,kindermann2020}.

\appendix
\setcounter{section}{0}  

\section{Appendix}\label{sec:appA}
\renewcommand{\theequation}{\thesection.\arabic{equation}} 
\makeatletter  
\@addtoreset{equation}{section} 
\makeatother

\newtheorem{appendixtheorem}{Theorem}[section]  
\newtheorem{appendixproposition}[appendixtheorem]{Proposition}
\newtheorem{appendixlemma}[appendixtheorem]{Lemma}
\renewcommand{\theappendixtheorem}{\thesection.\arabic{appendixtheorem}}
\renewcommand{\theappendixproposition}{\thesection.\arabic{appendixtheorem}}
\renewcommand{\theappendixlemma}{\thesection.\arabic{appendixtheorem}}

This appendix discusses the technical details of the theory presented. Our analysis is carried
out in infinite-dimensional spaces. We mainly extend techniques and results from
\cite{neubauer1988a,king1988} to the case of iterated Tikhonov. To the 
best of our knowledge, this is the first exploration of such an extension. Since this 
generalization poses nontrivial challenges and holds potential benefits for the analysis 
of other Tikhonov-type iterative methods, we believe that it deserves a dedicated section
of its own. We use the same notation as in~\cite{neubauer1988a}. To enhance readability, we 
provide Table \ref{tabnot}, which connects the notation of this appendix with the 
notation of Section \ref{sec:iAT}.

\begin{table}[h!]
    \centering
    \begin{tabular}{|c|c|c|c|c|c|c|}
        \hline
        Notation in Appendix \ref{sec:appA} & \( T \) & \( T_h \) & \(  Q_m  \) & \(  T_{h,m}\coloneqq Q_mT_h \) & \(  x^{\dagger} \) & \(  x_{\alpha,m,i}^{\delta,h} \) \\
        \hline
        Notation in Sections \ref{sec:AT} and \ref{sec:iAT} & \( T_n \) & \( T_n^{(\ell)} \) & \( \mathcal{R}_{\ell} \) & \( T_n^{(\ell)} \) & \(  x_n^{\dagger} \) & \(  x_{\alpha,n,i}^{\delta,\ell} \) \\
        \hline
    \end{tabular}
    \caption{Comparison of notations}\label{tabnot}
\end{table}

\subsection{Error estimates}\label{ssec:error_estimates}
Denote by $\mathcal{L}(\mathcal{X},\mathcal{Y})$ the space of linear operators from $\mathcal{X}$ to $\mathcal{Y}$. Let $T\in\mathcal{L}({\mathcal{X}},\mathcal{\mathcal{Y}})$ be a bounded linear operator and let $T_{h}$ be an approximation of $T$. We define, using the iterated Tikhonov~(iT) 
method, the computed approximate solution
\begin{equation}\label{it}\tag{iT}
x_{\alpha,i}\coloneqq\sum\limits_{k=1}^i\alpha^{k-1}(T^{\ast}T+\alpha I)^{-k}T^{\ast}y
\end{equation}
of~\eqref{eq}. When using $T_{h}$ instead of $T$, we obtain similarly $x_{\alpha,i}^{h}$.

\begin{appendixlemma}\label{lemma1}
Let $x_{\alpha,i}$ be defined as in~\eqref{it}. Then 
\begin{equation}\label{ineq1}
\|x_{\alpha,i}-x_{\alpha,i}^{h}\|\leq i(i+1)\frac{\|T-T_{h}\|}
{2\sqrt{\alpha}}\|x^{\dagger}\|+i\frac{\|Q_{h}(I-Q)y\|}{2\sqrt{\alpha}},
\end{equation}
where $Q$ and $Q_{h}$ are the orthogonal projectors on $\overline{\Range(T)}$ and 
$\overline{\Range(T_{h})}$, respectively.
\end{appendixlemma}

\begin{proof}
We have
\small
\begin{equation}\label{eq:difference}
  x_{\alpha,i}-x_{\alpha,i}^{h}=\sum\limits_{k=1}^{i}\alpha^{k-1}((T^{\ast}T+\alpha I)^{-k}-(T_h^{\ast}T_h+\alpha I)^{-k})T^{\ast}y+\sum\limits_{k=1}^i\alpha^{k-1}(T_h^{\ast}T_h+\alpha I)^{-k}(T^{\ast}-T_{h}^{\ast})y.
\end{equation}
\normalsize
Using the simple algebraic identities for $A,B\in\mathcal{L}(\mathcal{X},\mathcal{Y})$, 
\begin{equation*}
    A^{-k}-B^{-k}=B^{-k}(B^{k}-A^{k})A^{-k},\hspace{1cm}B^{k}-A^{k}=\sum\limits_{j=0}^{k-1}B^{j}(B-A)A^{k-1-j},
\end{equation*}
the first term in the right-hand side of~\eqref{eq:difference} can be written as
\begin{equation*}
\sum_{k=1}^{i}\sum\limits_{j=0}^{k-1}\alpha^{k-1}(T_h^{\ast}T_h+
\alpha I)^{j-k}(T_h^{\ast}T_h-T^{\ast}T)(T^{\ast}T+\alpha I)^{-1-j}T^{\ast}Qy,
\end{equation*}
and using the fact that 
$T_h^{\ast}T_h -T^{\ast}T=T_h^{\ast}(T_h-T)+(T_h^{\ast}-T^{\ast})T$,
we can split the above sum into
\begin{subequations}
\begin{equation}
\sum_{k=1}^{i}\sum\limits_{j=0}^{k-1}\alpha^{k-1}(T_h^{\ast}T_h+
\alpha I)^{j-k}T_h^{\ast}(T_h-T)(T^{\ast}T+\alpha I)^{-1-j}T^{\ast}Qy\label{term1}
\end{equation}
\begin{equation}
+\sum_{k=1}^{i}\sum\limits_{j=0}^{k-1}\alpha^{k-1}(T_h^{\ast}T_h+
\alpha I)^{j-k}(T_h^{\ast}-T^{\ast})T(T^{\ast}T+\alpha I)^{-1-j}T^{\ast}Qy.\label{term1b}
\end{equation}
\end{subequations}
Collecting factors $\alpha^{k-1}(T_h^{\ast}T_h+\alpha I)^{-k}$ from~\eqref{term1b} and the second 
term of~\eqref{eq:difference}, and using that $Qy=Tx^{\dagger}$ and 
$T^{\ast}_{h}=T^{\ast}_{h}Q_{h}$, their sum can be written as
\begin{subequations}
\begin{equation}
\sum_{k=1}^{i}\alpha^{k-1}(T_h^{\ast}T_h+\alpha I)^{-k}
\biggl(\sum_{j=0}^{k-1}(T_h^{\ast}T_h+\alpha I)^{j}(T_h^{\ast}-T^{\ast})
T(T^{\ast}T+\alpha I)^{-1-j}T^{\ast}-(T_{h}^{\ast}-T^{\ast})\biggr)Tx^{\dagger}
\label{term4a}
\end{equation}
\begin{equation}
+\sum_{k=1}^i\alpha^{k-1}(T_h^{\ast}T_h+\alpha I)^{-k}T_{h}^{\ast}Q_h(Q-I)y\label{term4}.
\end{equation}
\end{subequations}
Now rewrite the term~\eqref{term4a} as
\begin{subequations}
\begin{equation}
\sum_{k=1}^{i}\alpha^{k-1}(T_h^{\ast}T_h+\alpha I)^{-k}
\biggl(-\alpha(T_{h}^{\ast}-T^{\ast})(TT^{\ast}+\alpha I)^{-1}\biggr)Tx^{\dagger}
\label{term2}
\end{equation}
\begin{equation}
+\sum\limits_{k=1}^{i}\alpha^{k-1}(T_h^{\ast}T_h+\alpha I)^{-k}
\biggl(\sum\limits_{j=1}^{k-1}(T_h^{\ast}T_h+\alpha I)^{j}(T_h^{\ast}-T^{\ast})
T(T^{\ast}T+\alpha I)^{-1-j}T^{\ast}\biggr)Tx^{\dagger}\label{term3}.
\end{equation}
\end{subequations}
We obtain 
\begin{equation*}
    \|x_{\alpha,i}-x_{\alpha,i}^{h}\|\leq\|\eqref{term1}\|+\|\eqref{term2}\|+\|\eqref{term3}\|+\|\eqref{term4}\|.
\end{equation*}

Collecting $\|T-T_{h}\|\|x^{\dagger}\|$ from the first three terms and $\|Q_{h}(Q-I)y\|$ 
from the last term, and using the bounds
\begin{equation*}
\|(A^{\ast}A+\alpha I)^{-1}A^{\ast}A\|\leq 1,\qquad
\|(A^{\ast}A+\alpha I)^{-1}A^{\ast}\|\leq 1/(2\sqrt{\alpha})
\end{equation*}
and, similarly by switching the position of $\ast$, inequality \eqref{ineq1} follows.

For example, consider the argument of the nested sum~\eqref{term1}: Since $k-j-1\geq 0$ 
and $j\geq 0$, its norm is bounded as 
\small
\begin{align*}
 &\alpha^{k-1}\|(T_h^{\ast}T_h+\alpha I)^{j-k}T_h^{\ast}(T_h-T)(T^{\ast}T+
 \alpha I)^{-1-j}T^{\ast}Qy\|\\
 &\leq\alpha^{k-1}\|(T_h^{\ast}T_h+\alpha I)^{j-k+1}\|\|(T_h^{\ast}T_h+
 \alpha I)^{-1}T_h^{\ast}\|\|T-T_h\|\|(T^{\ast}T+\alpha I)^{-j}\|
 \|(T^{\ast}T+\alpha I)^{-1}T^{\ast}T\|\|x^{\dagger}\|\\
 &\leq\alpha^{k-1}\alpha^{j-k+1}(1/(2\sqrt{\alpha}))\|T-T_h\|\alpha^{-j}\|x^{\dagger}\|=
 \|T-T_h\|\|x^{\dagger}\|/(2\sqrt{\alpha}).
\end{align*}
\normalsize
Adding up all the terms, it follows that 
\begin{equation*}
\|\eqref{term1}\|+\|\eqref{term2}\|+\|\eqref{term3}\| \leq 
\frac{\|T-T_h\|}{2\sqrt{\alpha}}\|x^{\dagger}\|2\sum_{k=1}^ik =
i(i+1)\frac{\|T-T_h\|}{2\sqrt{\alpha}}\|x^{\dagger}\|.
\end{equation*}
\end{proof}

Consider a family $\lbrace W_m\rbrace_{m\in\N}$ of finite-dimensional subspaces of 
$\mathcal{Y}$ such that the orthogonal projector $Q_m$ into $W_m$ converges to $I$ on 
$\overline{\Range(T)}$. Define $x_{\alpha,m,i}^{h}$ and $x_{\alpha,m,i}^{\delta,h}$ 
similarly as in~\eqref{it}, i.e., by replacing $T_h$ by $T_{h,m}\coloneqq Q_mT_h$ and using
$y$ and $y^{\delta}$ as inputs, respectively. We now recall the useful 
\cite[Assumptions 2.3]{neubauer1988a}:

\begin{ass}
Let $W_m\subset\mathcal{Y}$, $T_h\in \mathcal{L}(\mathcal{X},\mathcal{Y})$, $\delta\geq 0$, and 
$y^{\delta}\in\mathcal{Y}$ be such that
\begin{equation}\label{conditions}\tag{C}
\begin{split}
    W_m\subset\overline{\Range(T)},\qquad &\Range(Q_mT_h)=W_m,\\
    \|Q_m(T-T_h)\|\leq h,\qquad &\|Q_m(y-y^{\delta})\|\leq\delta.
\end{split}
\end{equation}
\end{ass}

The next two lemmas are preparatory to the proof of Proposition~\ref{conv}.
\begin{appendixlemma}\label{lembound}
Under Conditions~\ref{conditions}, we have
\begin{equation*}
\|x^{\dagger}-x_{\alpha,m,i}^{\delta,h}\|\leq\alpha^i
\|(T^{\ast}_mT_m+\alpha I)^{-i}x^{\dagger}\|+(i(i+1)h\|x^{\dagger}\|+
i\delta)/(2\sqrt{\alpha}),
\end{equation*}
where $T_m\coloneqq Q_mT$.
\end{appendixlemma}
\begin{proof}
Let $x_{\alpha,m,i}$ be defined as in~\eqref{it} using $T_m\coloneqq Q_mT$. Then
\begin{equation*}
    \|x^{\dagger}-x_{\alpha,m,i}^{\delta,h}\|\leq\|x^{\dagger}-x_{\alpha,m,i}\|+\|x_{\alpha,m,i}-x_{\alpha,m,i}^{h}\|+\|x_{\alpha,m,i}^{h}-x_{\alpha,m,i}^{\delta,h}\|.
\end{equation*}
For the last term, we have the upper bound
\begin{equation*}
    \|x_{\alpha,m,i}^{h}-x_{\alpha,m,i}^{\delta,h}\|\leq\sum_{k=1}^i\|(T_{h,m}^{\ast}T_{h,m}+\alpha I)^{-1}T_{h,m}^{\ast}\|\|y-y^{\delta}\|\leq i\delta/(2\sqrt{\alpha}).
\end{equation*}
For the first term, using $T_mx^{\dagger}=Q_my$, we have
\begin{align}\label{i}
\|x^{\dagger}-x_{\alpha,m,i}\|&=\Bigl\|x^{\dagger}-\sum_{k=1}^i\alpha^{k-1}(T_m^{\ast}T_m+
\alpha I)^{-k}T_m^{\ast}T_mx^{\dagger}\Bigr\|\notag\\
 &=\Bigl\|\alpha(T_m^{\ast}T_m+\alpha I)^{-1}x^{\dagger}-
 \sum_{k=2}^i\alpha^{k-1}(T_m^{\ast}T_m+\alpha I)^{-k}T_m^{\ast}T_mx^{\dagger}\Bigr\|\notag\\
 &= \ldots\notag\\
 &=\alpha^{i}\|(T_m^{\ast}T_m+\alpha I)^{-i}x^{\dagger}\|.
\end{align}
The upper bound for the second term is derived using Lemma~\ref{lemma1} by replacing $T$ by
$T_m$, and $T_h$ by $T_{h,m}$. In accordance with Conditions~\ref{conditions}, this 
yields
\begin{equation*}
    \|x_{\alpha,m,i}-x_{\alpha,m,i}^{h}\|\leq i(i+1)\frac{h}{2\sqrt{\alpha}}\|x^{\dagger}\|.
\end{equation*}
Combining the three bounds we have computed thus far, the lemma follows.
\end{proof}

For $\nu\geq 0$ and $\rho>0$, we introduce the set 
\begin{equation*}
    \mathcal{X}_{\nu,\rho}\coloneqq\lbrace x\in\mathcal{X} \mid x=(T^{\ast}T)^{\nu}w,\; 
w\in \ker(T)^\perp \mbox{ and } \|w\|\leq\rho\rbrace.
\end{equation*}
Assuming that $x^{\dagger}\in\mathcal{X}_{\nu,\rho}$, we can bound~\eqref{i} by 
\begin{equation*}
    \|x^{\dagger}-x_{\alpha,m,i}\|\leq\alpha^{i}\|(T_m^{\ast}T_m+\alpha I)^{-i}(T_m^{\ast}T_m)^{\nu}w\|+\|[(T^{\ast}T)^{\nu}-(T_m^{\ast}T_m)^{\nu}]w\|.
\end{equation*}
Define 
\begin{equation*}
    \gamma_m\coloneqq\|(I-Q_m)T\|.
\end{equation*}

\begin{appendixlemma}\label{lemlim}
Let $x^{\dagger}\in\mathcal{X}_{\nu,\rho}$ and let $w$ be fixed. Defining
\begin{equation*}
 b(\alpha,m)\coloneqq \alpha^{i(1-\nu)}\|(T_m^{\ast}T_m+\alpha I)^{-i}
 (T_m^{\ast}T_m)^{\nu}w\|,
\end{equation*}
it holds that 
\begin{equation*}
    \lim_{\substack{m\rightarrow\infty \\ \alpha\rightarrow 0^{+}}}b(\alpha,m)=0
\end{equation*}
if $\nu\in [0,1)$ and 
\begin{equation}\label{condplus}
\gamma_m=o\Bigl(\mu_{m}^{\frac{\nu(i-1)}{1-\nu}}\Bigr),
\end{equation}
for $\mu_{m}\coloneqq\min(\sigma(T^{\ast}_mT_m)\setminus\lbrace 0\rbrace)$, where $\sigma(T_m^{\ast}T_m)$ denotes the spectrum of $T_m^{\ast}T_m$.
\end{appendixlemma}

\begin{proof}
Letting $\lbrace E_{\mu}^{m}\rbrace_{\mu\in\R}$ be a spectral family for $T^{\ast}_mT_m$ (see e.g.~\cite[Section 2.3]{engl1996}) 
and $f_{\nu}(\mu)=\alpha^i(\mu+\alpha)^{-i}\mu^{\nu}$, we have the equality 
\begin{equation*}
  b^{2}(\alpha,m)=\int_0^{\infty}\alpha^{-2i\nu}f_{\nu}^{2}(\mu)dE_{\mu}^{m}\|w\|^2.
\end{equation*}
By H\"older's inequality,
\begin{align*}
 b^2(\alpha,m)&\leq\left[\int_0^{\infty}\left(\frac{\alpha}{\mu+\alpha}\right)^{2i}
dE_{\mu}^{m}\|w\|^2\right]^{1-\nu}\left[\int_0^{\infty}
\left(\frac{\mu}{(\mu+\alpha)^{i}}\right)^{2}dE_{\mu}^{m}\|w\|^2\right]^{\nu}\\
 &\leq\left[\int_0^{\infty}\left(\frac{\alpha}{\mu+\alpha}\right)^{2i}dE_{\mu}^{m}
 \|w\|^2\right]^{1-\nu}\mu_{m}^{2\nu(1-i)}\|w\|^{2\nu}.
\end{align*}
Now
\begin{equation*}
 b^{\frac{1}{1-\nu}}(\alpha,m) \leq\mu_{m}^{\frac{\nu(1-i)}{1-\nu}}\alpha\|(T_m^{\ast}T_m+
 \alpha I)^{-1}w\|\|w\|^{\frac{\nu}{1-\nu}}
 \stackrel{\substack{\alpha\to 0 \\ m\to\infty}}{\longrightarrow} 0
\end{equation*}
follows from the proof of~\cite[Lemma 2.4]{king1988}.
\end{proof}

\begin{appendixproposition}\label{conv}
Let Conditions~\ref{conditions} and~\eqref{condplus} be satisfied. Moreover, let 
$x^{\dagger}\in\mathcal{X}_{\nu,\rho}$. For $m$ sufficiently large, there is a unique
$\alpha>0$ that solves
\begin{equation*}
 \alpha^{\frac{2i+1}{2}}\|(T_m^{\ast}T_m+\alpha I)^{-i}(T_m^{\ast}T_m)^{\nu}w\|=
 (i(i+1)h\|x^{\dagger}\|+i\delta)/2.
\end{equation*}
Moreover, it holds $\alpha\rightarrow 0$ for $m\rightarrow\infty$, $h,\delta\rightarrow 0$, and if
$\sqrt{\alpha}=o\Bigl(\mu_m^{\frac{\nu(i-1)}{1-\nu}}\Bigr)$, then we have
\begin{equation*}
  \|x^{\dagger}-x_{\alpha,m,i}^{\delta,h}\|=\begin{cases}
      o(1)&\text{if $\nu=0$,}\\
      o((h+\delta)^{\frac{2i\nu}{2i\nu+1}})+O(\gamma_m^{2\nu}\|w\|)&\text{if $0<\nu<1$,}\\
      O((h+\delta)^{\frac{2i}{2i+1}})+O(\gamma_m\|(I-Q_m)Tw\|)&\text{if $\nu=1$.}
  \end{cases}
\end{equation*}
\end{appendixproposition}

\begin{proof}
Follows from Lemma~\ref{lembound} and Lemma~\ref{lemlim} analogously to~\cite[Proposition 2.6]{neubauer1988a} and~\cite[Theorem 2.1]{king1988}.
\end{proof}


\subsection{The parameter choice method}\label{parasele}
We describe the a posteriori parameter choice method that is used in Algorithm~\ref{algo:itarntik} 
in Section~\ref{sec:iAT}. Define
\begin{equation*}
 f(m,\alpha,w,A)\coloneqq\alpha^{2i+1}\langle(A_mA^{\ast}_m+
 \alpha I)^{-2i-1}Q_mw,Q_mw\rangle,
\end{equation*}
where $A_m\coloneqq Q_mA$ and $\langle\cdot,\cdot\rangle$ denotes the Euclidean scalar product.

\begin{appendixproposition}
Let $C$ and $E$ be positive constants such that
\begin{equation}\label{condEC}
  Eh+C\delta\leq\|Q_my^{\delta}\|.
\end{equation}
Then there is a unique $\alpha>0$ that solves 
\begin{equation}\label{select}
  f(m,\alpha,y^{\delta},T_h)=(Eh+C\delta)^2.
\end{equation}
\end{appendixproposition}

\begin{proof}
Analogous to~\cite[Proposition 3.1]{neubauer1988a}.
\end{proof}

We now extend the result of~\cite[Lemma 3.2]{neubauer1988a}, which will be used in Theorem \ref{conv2}.

\begin{appendixlemma}\label{lemAB}
Let $A,B\in \mathcal{L}(\mathcal{X},\mathcal{Y})$. Then
\begin{equation*}
 \alpha^{\frac{2i+1}{2}}\|[(AA^{\ast}+\alpha I)^{-\frac{2i+1}{2}}-
 (BB^{\ast}+\alpha I)^{-\frac{2i+1}{2}}]A\|\leq (2i+1)\|A-B\|.
\end{equation*}
\end{appendixlemma}

\begin{proof}
We rewrite the element in the square brackets as
\begin{align*}
 &\sum\limits_{k=0}^{i-1}(BB^{\ast}+\alpha I)^{k-\frac{2i+1}{2}}(BB^{\ast}-
  AA^{\ast})(AA^{\ast}+\alpha I)^{-1-k}\\
 &+(BB^{\ast}+\alpha I)^{-\frac{1}{2}}((BB^{\ast}+\alpha I)^{\frac{1}{2}}-
 (AA^{\ast}+\alpha I)^{\frac{1}{2}})(AA^{\ast}+\alpha I)^{-\frac{2i+1}{2}}.
\end{align*}
For the first term, we use $BB^{\ast}-AA^{\ast}=B(B^{\ast}-A^{\ast})+(B-A)A^{\ast}$, obtaining
\begin{equation*}
 \alpha^{\frac{2i+1}{2}}\Bigl\|\sum_{k=0}^{i-1}(BB^{\ast}+\alpha I)^{k-\frac{2i+1}{2}}
 (BB^{\ast}-AA^{\ast})(AA^{\ast}+\alpha I)^{-1-k}A\Bigr\|\leq 2i\|A-B\|.
\end{equation*}
The second term is bounded from above by
\begin{equation*}
 \alpha^{\frac{3}{2}}\|(BB^{\ast}+\alpha I)^{-\frac{1}{2}}
 ((BB^{\ast}+\alpha I)^{\frac{1}{2}}-(AA^{\ast}+\alpha I)^{\frac{1}{2}})
 (AA^{\ast}+\alpha I)^{-\frac{3}{2}}A\|
\end{equation*}
and the lemma follows from the proof of~\cite[Lemma 3.2]{neubauer1988a}.
\end{proof}

\begin{appendixtheorem}\label{conv2}
Let $C>1$ and $E>(2i+1)\|x^{\dagger}\|$. For each quintuple 
$(W_m,h,T_h,\delta,y^{\delta})$, let
Conditions~\ref{conditions},~\eqref{condplus}, and~\eqref{condEC} hold. 
Moreover, let $x^{\dagger}\in\mathcal{X}_{\nu,\rho}$ and let $\alpha>0$ be the unique 
solution of~\eqref{select}. Then, if 
$\sqrt{\alpha}=o\Bigl(\mu_m^{\frac{\nu(i-1)}{1-\nu}}\Bigr)$, the same asymptotic estimates 
of Proposition~\ref{conv} hold for this parameter choice.
\end{appendixtheorem}

\begin{proof}
The result follows from Lemma~\ref{lemAB} analogously to~\cite[Lemma 3.3, Proposition 3.4, Theorem 3.5]{neubauer1988a}.
\end{proof}

We now extend the result of~\cite[Proposition 3.6]{neubauer1988a}.

\begin{appendixproposition}\label{comp}
Let $C=1$ and $E=\|x^{\dagger}\|$. Assume that 
\eqref{condEC} holds and let $\alpha>0$ be the unique solution of~\eqref{select}. 
Then, for all $\tilde{\alpha}\geq\alpha$, it holds that 
$\|x^{\dagger}-x_{\alpha,m,i}^{\delta,h}\|\leq
\|x^{\dagger}-x_{\tilde{\alpha},i,m}^{\delta,h}\|$.
\end{appendixproposition}

\begin{proof}
Let $\lbrace F_{\mu}^{h,m}\rbrace_{\mu\in\R}$ be a spectral family for 
$T_{h,m}T^{\ast}_{h,m}$ and define 
$e(\alpha)\coloneqq\frac{1}{2}\|x^{\dagger}-x_{\alpha,m,i}^{\delta,h}\|^2$. Then
\begin{equation*}
  \frac{de(\alpha)}{d\alpha}=i\biggl\langle T_{h,m}x^{\dagger}-
  \int_0^{\infty}\frac{(\mu+\alpha)^i-\alpha^i}{(\mu+\alpha)^i}dF_{\mu}^{h,m}
  Q_my^{\delta},\int_{0}^{\infty}\frac{\alpha^{i-1}}{(\mu+\alpha)^{i+1}}dF_{\mu}^{h,m}
  Q_my^{\delta}\biggr\rangle.
\end{equation*}
Adding and subtracting $i\alpha^{2i-1}\bigl\|(T_{h,m}T_{h,m}^{\ast}+
\alpha I)^{\frac{-2i-1}{2}}Q_my^{\delta}\bigr\|^2$, we obtain
\begin{align*}
 \frac{de(\alpha)}{d\alpha}&=i\alpha^{2i-1}\bigl\|(T_{h,m}T_{h,m}^{\ast}+\alpha I)^{\frac{-2i-1}{2}}
 Q_my^{\delta}\bigr\|^2\\
 &+i\biggl\langle\int_{0}^{\infty}\frac{\alpha^{i-1}}{(\mu+\alpha)^{\frac{1}{2}}}
 dF_{\mu}^{h,m}(T_{h,m}x^{\dagger}-Q_my^{\delta}),(T_{h,m}T_{h,m}^{\ast}+
 \alpha I)^{\frac{-2i-1}{2}}Q_my^{\delta}\biggr\rangle.
\end{align*}
Thus, collecting $i\bigl\|(T_{h,m}T_{h,m}^{\ast}+\alpha I)^{\frac{-2i-1}{2}}
Q_my^{\delta}\bigr\|\eqqcolon K$ from the two terms above, we have
\begin{equation*}
 \frac{de(\alpha)}{d\alpha}\geq K\left(\alpha^{2i-1}\bigl\|(T_{h,m}T_{h,m}^{\ast}+\alpha I)^{\frac{-2i-1}{2}}
 Q_my^{\delta}\bigr\|-\biggl\|\int_{0}^{\infty}\frac{\alpha^{i-1}}
 {(\mu+\alpha)^{\frac{1}{2}}}dF_{\mu}^{h,m}(T_{h,m}x^{\dagger}-Q_my^{\delta})\biggr\|\right).
\end{equation*}
The proposition now follows from
\begin{equation*}
  \biggl\|\int_{0}^{\infty}\frac{\alpha^{i-1}}{(\mu+\alpha)^{\frac{1}{2}}}dF_{\mu}^{h,m}
  (T_{h,m}x^{\dagger}-Q_my^{\delta})\biggr\|\leq\alpha^{i-\frac{3}{2}}
  \|T_{h,m}x^{\dagger}-Q_my^{\delta}\|
\end{equation*}
similarly as~\cite[Proposition 3.6]{neubauer1988a}.
\end{proof}

Analogously to~\cite[Remark 3.7]{neubauer1988a}, if an estimation of $\|x^{\dagger}\|$ is 
not available, then the function $D\|x_{\alpha,m,i}^{\delta,h}\|$ with a constant 
$D\geq 1$ can be used in place of $\|x^{\dagger}\|$ to define the constant $E$. It can be shown 
that for $\alpha$ satisfying~\eqref{select}, this choice yields the same rates 
of convergence as presented in Theorem~\ref{conv2}.


\subsection{Other convergence rates}\label{sec:appB}
We derive results on the convergence rates similar to the ones of 
\cite{yxps2022,gfrerer1987,scherzer1993} for~\eqref{it}, and provide a brief description  
of all cases $\nu\geq 0$.
The main difference here is that we no longer require the hypothesis \eqref{condplus} on $\mu_m$, i.e., on the spectrum of $T^{\ast}_mT_m$, to be satisfied, which leads to slower convergence rates.

\begin{appendixlemma}\label{lemlim2}
If $\nu\in [0,i)$, then it holds
\begin{equation*}
  \lim\limits_{\substack{m\rightarrow\infty \\ \alpha\rightarrow 0^{+}}}\alpha^{i-\nu}\|(T_m^{\ast}T_m+\alpha I)^{-i}(T_m^{\ast}T_m)^{\nu}w\|=0.
\end{equation*}
\end{appendixlemma}
\begin{proof}
Analogous to Lemma~\ref{lemlim}.
\end{proof}

\begin{appendixproposition}\label{convclass}
Let Conditions~\ref{conditions} be satisfied and let 
$x^{\dagger}\in\mathcal{X}_{\nu,\rho}$. Let $\alpha>0$ be as in Proposition~\ref{conv}. 
Then we have
\begin{equation*}
  \|x^{\dagger}-x_{\alpha,m,i}^{\delta,h}\|=\begin{cases}
      o(1)&\text{if $\nu=0$,}\\
      o((h+\delta)^{\frac{2\nu}{2\nu+1}})+O(\gamma_m^{2\nu}\|w\|)&\text{if $0<\nu<i$,}\\
      O((h+\delta)^{\frac{2i}{2i+1}})+O(\gamma_m\|(I-Q_m)Tw\|)&\text{if $\nu=i$.}
  \end{cases}
\end{equation*}
\end{appendixproposition}

\begin{proof}
Follows from Lemma~\ref{lemlim2} analogously to Proposition~\ref{conv}.
\end{proof}

Results similar to those of Section~\ref{parasele} on the parameter choice method, and of 
Proposition~\ref{convf} and its corollaries, remain valid for $\nu\in[0,i]$. For the case
$\nu>i$, we have the following result.

\begin{appendixproposition}\label{bignu}
Let Conditions~\ref{conditions} hold and let $x^{\dagger}\in\mathcal{X}_{\nu,\rho}$
with $\nu>i$. For $\alpha=(h+\delta)^{\frac{2}{2i+1}}$, we have
\begin{equation*}
  \|x^{\dagger}-x_{\alpha,m,i}^{\delta,h}\|=O((h+\delta)^{\frac{2i}{2i+1}})+
  \|[(T^{\ast}T)^{\nu}-(T_m^{\ast}T_m)^{\nu}]w\|.
\end{equation*}
\end{appendixproposition}

\begin{proof}
Analogous to the first part of Proposition~\ref{conv}.
\end{proof}

\section*{Acknowledgments}
We would like to thank a referee for comments that lead to clarifications of the
presentation.
The work of the first author is partially supported by NSFC (Grant No.~12250410253) and by the Startup Fund of Sun~Yat-sen~University. The work of the second author is partially supported by MIUR - PRIN 2022 N.2022ANC8HL and GNCS-INdAM.


\bibliography{Article_rev}

\end{document}